\documentclass[twoside,11pt]{article}

\usepackage[utf8]{inputenc}
\usepackage[T1]{fontenc}
\usepackage[table]{xcolor}
\usepackage[abbrvbib,preprint]{jmlr2e}

\usepackage{lastpage}
\usepackage{graphicx}
\usepackage{amsmath,amsfonts,amsopn,mathrsfs}
\usepackage{newtxtext,newtxmath}
\usepackage{tikz}
\usetikzlibrary{arrows.meta,positioning,fit,calc}
\usepackage{array}
\usepackage{enumitem}
\usepackage{bm}
\usepackage{url}
\usepackage{booktabs}
\usepackage{multirow}
\usepackage{nicefrac}
\usepackage{microtype}
\usepackage{algorithm}
\usepackage{algorithmic}
\usepackage[misc,geometry]{ifsym}
\usepackage{etoolbox}
\usepackage{makecell}
\usepackage{tabularx}
\allowdisplaybreaks
\usepackage{threeparttable}
\usepackage{placeins}

\newcommand{\E}{\mathbb{E}}
\newcommand{\R}{\mathbb{R}}
\newcommand{\ind}{\mathbf{1}}
\DeclareMathOperator{\Pro}{\mathbb{P}}

\newtheorem{thm}{Theorem}
\newtheorem{lem}{Lemma}

\newtheorem{pros}{Proposition}

\newtheorem{rem}{Remark}

\newtheorem{assumpt}{Assumption}
\newenvironment{namedproof}[1]{\par\noindent{\bf #1\ }}{\hfill\BlackBox\\[2mm]}

\usepackage{placeins}

\usepackage[nameinlink,noabbrev,capitalize]{cleveref}
\hypersetup{
  hypertexnames=false,
  pdftitle={Long-Run Dynamics of AdaGrad-Norm: Trajectory Stability and Asymptotic Stationarity},
  pdfauthor={Ruinan Jin and Xiaoyu Wang},
  pdfsubject={Long-run dynamics, trajectory stability, and asymptotic stationarity of AdaGrad-Norm},
  pdfkeywords={AdaGrad-Norm, nonconvex optimization, stochastic approximation, trajectory stability, adaptive stepsizes}
}
\crefname{thm}{Theorem}{Theorems}
\crefname{lem}{Lemma}{Lemmas}
\crefname{prop}{Property}{Properties}
\crefname{cor}{Corollary}{Corollaries}
\crefname{assumpt}{Assumption}{Assumptions}
\crefname{rem}{Remark}{Remarks}
\crefname{exa}{Example}{Examples}
\crefname{pros}{Proposition}{Propositions}

\jmlrheading{1}{2026}{1-\pageref{LastPage}}{--}{--}{--}{Ruinan Jin and Xiaoyu Wang}
\ShortHeadings{Long-Run Dynamics of AdaGrad-Norm}{Jin and Wang}
\firstpageno{1}

\begin{document}

\title{Long-Run Dynamics of AdaGrad-Norm: Trajectory Stability and Asymptotic Stationarity}

\author{\name Ruinan Jin \email jinruinan3@gmail.com \\
       \addr Department of Electrical and Computer Engineering\\
       The Ohio State University\\
       Columbus, OH 43210, USA
       \AND
       \name Xiaoyu Wang\textsuperscript{*} \email wangxiaoyu@ucas.ac.cn \\
       \addr School of Mathematical Sciences\\
       University of Chinese Academy of Sciences\\
       No.~19A Yuquan Road, Shijingshan District\\
       Beijing 100049, China}

\editor{To be assigned}

\maketitle

\renewcommand{\thefootnote}{\fnsymbol{footnote}}
\setcounter{footnote}{1}
\footnotetext{Corresponding author: Xiaoyu Wang, School of Mathematical Sciences, University of Chinese Academy of Sciences, No.~19A Yuquan Road, Shijingshan District, Beijing 100049, China. Email: \texttt{wangxiaoyu@ucas.ac.cn}.}
\setcounter{footnote}{0}
\begin{abstract}%
AdaGrad-Norm is widely studied, yet its long-run behavior in smooth nonconvex
optimization remains delicate. For the standard current-denominator recursion,
the same stochastic gradient determines both the update and its normalizer,
rendering the effective stepsize nonpredictable. At the square-root
normalization, the quadratic smoothness budget admits only logarithmic control
rather than a uniform summability bound. We establish trajectory stability and
full-sequence asymptotic stationarity directly for this unmodified recursion.
Under global smoothness, non-flatness at infinity, conditional unbiasedness,
and an affine conditional second-moment bound, we prove
\(\E[\sup_{n\geq1} g(\theta_n)]<\infty\) and
\(\E[\sup_{n\geq1}\|\theta_n\|]<\infty\), without assuming bounded
iterates. For stationarity, we use a local conditional
relative-moment condition that permits unbounded oracle values and covers
fixed-size finite-sum mini-batches and locally nondegenerate additive noise.
Together with a weak Sard condition, it yields
\(\|\nabla g(\theta_n)\|\to0\) almost surely and
\(\E\|\nabla g(\theta_n)\|^2\to0\).
The proof combines a current-sample Lyapunov correction, a stopped-excursion
stability argument, and finite-window control for the actual nonpredictable
stepsize. We also analyze power normalization with \(q>1/2\), where the
quadratic normalization budget is summable, identifying \(q=1/2\) as the
boundary of the proof mechanism developed here rather than a universal
algorithmic phase transition.
\end{abstract}


\section{Introduction}
\label{sec:introduction}

Adaptive gradient methods couple optimization dynamics with data-dependent
normalization. This coupling is particularly delicate for AdaGrad-Norm,
where a single scalar accumulator is updated from stochastic gradients and
simultaneously determines the learning rate. In the standard
current-denominator recursion, the same stochastic gradient both updates the
accumulator and moves the iterate. Understanding the resulting long-run
dynamics therefore requires controlling the interaction among stochastic
descent, self-normalization, and the evolution of the effective stepsize.

We study the standard, unmodified current-denominator AdaGrad-Norm recursion
\citep{ward2020adagrad},
\begin{equation}
\label{eq:adagrad-norm-intro}
 S_n=S_{n-1}+\|G_n\|^2,
 \qquad
 \theta_{n+1}=\theta_n-\frac{\alpha_0}{\sqrt{S_n}}G_n,
 \quad G_n:=\nabla g(\theta_n,\xi_n),
\end{equation}
with deterministic \(S_0,\alpha_0>0\).  Two features obstruct a direct
application of classical stochastic-approximation arguments.  First,
\(\alpha_0/\sqrt{S_n}\) depends on the current oracle output and is therefore
not measurable with respect to the past filtration.  Even under conditional
unbiasedness, the update cannot be decomposed directly into a predictable
drift and a martingale difference.  Second, smoothness generates the
quadratic normalization budget
\begin{equation}
\label{eq:intro-log-energy}
    \sum_{n=1}^{T}\frac{\alpha_0^2\|G_n\|^2}{S_n}
    \leq
    \alpha_0^2\log\frac{S_T}{S_0},
\end{equation}
which is finite on every finite horizon but need not be summable over an
infinite trajectory.  Thus, at the square-root normalization, neither a predictable-stepsize decomposition nor a summable quadratic normalization budget is directly available.

Existing finite-horizon analyses quantify optimization progress at a
prescribed computational budget, but such guarantees do not by themselves
provide uniform-in-time control of a realized trajectory. In particular,
they do not imply boundedness of the iterates, integrability of trajectory
suprema, or convergence of the full gradient-norm sequence. These are the
trajectory-level questions addressed here.

To address these trajectory-level questions, we take a stability-first route.
We first establish \(L^1\) control of the running objective and iterate
suprema, without assuming bounded iterates, and only then invoke the limiting
gradient-flow dynamics.  Two features are central to this route: non-flatness
converts objective control into iterate control, and finite-window
perturbations are handled for the \emph{actual} current-dependent stepsize
rather than for a predictable surrogate.  This distinguishes the analysis
from predictable-denominator arguments for past-gradient AdaGrad
\citep{li2019convergence} and from the current-denominator framework of
\citet{jin2022convergence}.  The additional stochastic condition used for
stationarity is local: a conditional relative-moment comparison near the
stationary region, which still permits unbounded oracle values and covers
standard finite-sum and additive-noise models.

\subsection{Contributions}
\label{subsec:intro-contributions}

The main contributions are as follows.
\begin{enumerate}[label=\textnormal{(\arabic*)},leftmargin=*,itemsep=2pt]
\item
\textbf{Integrable trajectory stability.}
Under global smoothness, non-flatness at infinity, conditional unbiasedness,
and an affine conditional second-moment bound, we prove
\[
 \E[\sup_n g(\theta_n)]<\infty,
 \qquad \E[\sup_n\|\theta_n\|]<\infty.
\]
Thus iterate boundedness is derived from objective stability through the
linear-growth consequence of non-flatness rather than assumed a priori.

\item
\textbf{Full-sequence stationarity under a local relative-moment condition.}
For stationarity, we impose a conditional \((2+\rho)\)-to-2 moment bound only
near the stationary region. It permits unbounded oracle values, holds for
fixed-size finite-sum mini-batches and locally nondegenerate additive noise,
and rules out an adapted-oracle counterexample showing that an ordinary
uniform local higher-moment bound is insufficient in the general model.
Together with weak Sard, it yields
\(\|\nabla g(\theta_n)\|\to0\) almost surely and
\(\E\|\nabla g(\theta_n)\|^2\to0\); the finite-window estimate is proved for
the actual current-dependent stepsize.

\item
\textbf{A stability mechanism for the square-root boundary.}
We combine a current-sample Lyapunov correction with a two-threshold,
three-phase stopped-excursion argument. For normalization \(S_n^{-q}\), the
quadratic budget is summable for \(q>1/2\), logarithmic at \(q=1/2\), and may
grow polynomially for \(q<1/2\). We also analyze the summable regime
\(q>1/2\), thereby isolating the square-root case as the boundary of the proof
mechanism developed here.
\end{enumerate}

\subsection{Related Work}
\label{subsec:intro-related}

\paragraph{Finite-horizon theory for AdaGrad-Norm.}
AdaGrad was introduced for online and stochastic convex optimization by
\citet{duchi2011adaptive}. For smooth nonconvex optimization,
\citet{ward2020adagrad} established finite-horizon guarantees for
AdaGrad-Norm, with subsequent work treating stronger geometries, affine or
unbounded-gradient noise models, high-probability guarantees, generalized
smoothness, and heavy-tailed noise
\citep{xie2020linear,faw2022power,wang2023convergence,kavis2022high,
liu2023high,attia2023adagrad,faw2023beyond,hong2024revisiting,
chezhegov2025clipping}. These results are complementary to the
infinite-horizon trajectory questions considered here.

\paragraph{Asymptotic and current-denominator analyses.}
For generalized AdaGrad stepsizes based on past stochastic gradients,
\citet{li2019convergence} proved almost-sure convergence of the gradient norm
under a predictable denominator and stronger global regularity and noise
conditions. Related asymptotic results include deterministic AdaGrad in smooth
convex optimization \citep{traore2021sequential}, analyses of current-denominator
AdaGrad(Norm) in convex and quasar-convex settings \citep{liu2023convergence},
and AdaGrad-/RMSProp-type stochastic approximation with decreasing stepsizes
\citep{gadat2022asymptotic}. Most closely related, \citet{jin2022convergence}
proved that current-denominator AdaGrad approaches a connected stationary
component almost surely under independent oracle samples, a bounded stationary
set with finitely many connected components, a local separation condition,
and an affine second-moment bound. In contrast, our analysis first derives
\(L^1\) control of the running objective and iterate suprema and then establishes
finite-window perturbation control for the actual current-dependent stepsize,
leading to full-sequence gradient-norm convergence almost surely and in mean
square. Broader last-iterate and almost-sure frameworks for stochastic adaptive
methods are developed in \citet{liang2026lastiterate}.

\paragraph{Stochastic approximation and adaptive optimization.}
Our ODE argument uses classical martingale, Lyapunov, and
asymptotic-pseudotrajectory tools
\citep{robbins1951stochastic,robbins1971convergence,
bertsekas2000gradient,kushner2003stochastic,borkar2008stochastic,
benaim2006dynamics,borkar2000ode}. The current denominator prevents the usual
predictable-stepsize martingale decomposition, while the local
higher-to-second moment comparison used here is related to martingale limit
theory \citep{machkouri2007exact,fan2019exact}. Analyses of Adam, RMSProp,
and coordinatewise AdaGrad involve additional momentum or preconditioning
structure and require different techniques
\citep{barakat2021convergence,JMLR:v25:23-0576,DBLP:conf/icml/JinLY025,
jiang2025adagrad}.

\begin{table}[t]
\centering
\caption{Nearest asymptotic comparisons for stochastic AdaGrad-type methods.
Common smoothness and unbiasedness conditions are omitted for brevity.
The remaining assumption sets are not totally ordered; the table highlights
structural differences rather than dominance.}
\label{tab:literature-comparison}
\renewcommand{\arraystretch}{1.15}
\setlength{\tabcolsep}{3pt}
\footnotesize
\begin{tabularx}{\textwidth}{@{}
>{\raggedright\arraybackslash}p{0.13\textwidth}
>{\raggedright\arraybackslash}p{0.11\textwidth}
>{\raggedright\arraybackslash}X
>{\raggedright\arraybackslash}p{0.20\textwidth}
>{\raggedright\arraybackslash}p{0.22\textwidth}@{}}
\toprule
Work
& Denominator
&  Key conditions
& Stability
& Convergence
\\
\midrule

\citet{li2019convergence}
& Past-gradient
& Generalized exponent \(1/2+\epsilon\); bounded-support noise

& No trajectory-supremum guarantee
& Gradient norm converges to zero a.s.
\\

\citet{jin2022convergence}
& Current-gradient
& Bounded stationary set with finitely many connected components;
local separation; affine second moment
& A.s. boundedness follows from convergence
& Distance to a connected stationary component tends to zero a.s.
\\

\rowcolor{gray!10}
\textbf{This work}
& \textbf{Current-gradient}
& Non-flatness; affine second moment; local relative-moment condition and weak Sard for stationarity
& \textbf{\(L^1\) control of objective and iterate suprema}
& \textbf{Gradient norm converges a.s. and in mean square}
\\

\bottomrule
\end{tabularx}
\end{table}

The remainder of the paper is organized as follows.  \Cref{sec:prelim}
formalizes the model and assumptions, \cref{sec:main-results} states and
interprets the main theorem, and \cref{sec:main-analysis} develops its
stability and stationarity analysis.  \Cref{sec:scalar-extensions} then places
the square-root case within the broader power-normalized family, followed by
the numerical diagnostics in \cref{sec:numerical} and the concluding
discussion in \cref{sec:conclusion}.

\section{Problem Setup and Assumptions}
\label{sec:prelim}

This section fixes notation and separates the assumptions used for trajectory
stability from those needed only for asymptotic stationarity.  The distinction
will be reflected explicitly in the layered main theorem in
\cref{sec:main-results}.

We minimize a continuously differentiable objective \(g:\R^d\to\R\) with
\(g_\star:=\inf_{\theta\in\R^d}g(\theta)>-\infty\).
Since additive constants change neither gradients nor iterates, we normalize
\(g_\star=0\), so \(g\geq0\); the infimum need not be attained.
We write
\(\Theta_{\mathrm{stat}}
 :=\{\theta\in\R^d:\nabla g(\theta)=0\}\)
and denote its critical-value set by
\(g(\Theta_{\mathrm{stat}})\).
We write \(\E[\cdot]\) and \(\E[\cdot\mid\mathscr F_n]\) for expectation
and conditional expectation, and \(\ind_A\) for the indicator of an event
\(A\).

\subsection{Objective and Oracle Assumptions}
\label{subsec:stability-assumptions}

\begin{assumpt}[Smoothness and non-flatness at infinity]
\label{ass_g_poi}
The objective \(g\) satisfies the following conditions.
\begin{enumerate}[label=\textnormal{(\roman*)},leftmargin=*]

    \item
    \label{ass_g_poi:i2}
    \textbf{Global smoothness.}
    There exists \(L>0\) such that
    \[
        \|\nabla g(\theta)-\nabla g(\theta')\|
        \leq
        L\|\theta-\theta'\|
        \qquad
        \text{for all }\theta,\theta'\in\R^d.
    \]

    \item
    \label{ass_g_poi:i3}
    \textbf{Non-flatness at infinity.}
    \[
        \liminf_{\|\theta\|\to\infty}
        \|\nabla g(\theta)\|
        >
        0.
    \]

\end{enumerate}
\end{assumpt}

Under \cref{ass_g_poi}~\ref{ass_g_poi:i3}, every sufficiently small
gradient sublevel set
\(K_\delta:=\{\theta:\|\nabla g(\theta)\|\leq\delta\}\)
is compact.  We will also use the following stronger consequence.

\begin{lem}[Linear growth from non-flatness]
\label{lem:nonflat-linear-growth}
Suppose that \cref{ass_g_poi} holds under the preceding normalization, and
let \(\ell_\infty>0\) denote the liminf in
\cref{ass_g_poi}~\ref{ass_g_poi:i3}.
For every \(c\in(0,\ell_\infty)\), there exists \(R_c<\infty\) such
that
\begin{equation}
\label{eq:nonflat-linear-growth}
    g(\theta)
    \geq
    c\bigl(\|\theta\|-R_c\bigr)_+
    \qquad\text{for every }\theta\in\R^d.
\end{equation}
In particular, \(g\) is coercive and every objective sublevel set is
compact.
\end{lem}

Thus non-flatness is stronger than coercivity: for example,
\(g(x)=\log(1+x^2)\) is coercive but has \(g'(x)\to0\).
Conversely, it holds automatically for
\(g(\theta)=h(\theta)+(\lambda/2)\|\theta\|^2\) when
\(\lambda>0\) and \(\nabla h\) is uniformly bounded, covering the
regularized finite-sum model used in \cref{subsec:numerical-california}.

Global smoothness and nonnegativity also imply
\begin{equation}
\label{eq:smooth-gradient-objective-bound}
    \|\nabla g(\theta)\|^2
    \leq
    2L g(\theta)
    \qquad
    \text{for all }\theta\in\R^d.
\end{equation}
This is the standard descent-lemma bound evaluated at
\(\theta-L^{-1}\nabla g(\theta)\).

Let \((\Omega,\mathscr F,\Pro)\) be a probability space equipped with
a filtration \(\{\mathscr F_n\}_{n\geq0}\). At iteration \(n\), the
iterate \(\theta_n\) and the pre-update accumulator \(S_{n-1}\) are
\(\mathscr F_{n-1}\)-measurable. After observing the oracle sample
\(\xi_n\), the stochastic gradient
\(G_n:=\nabla g(\theta_n,\xi_n)\) is
\(\mathscr F_n\)-measurable.

\begin{assumpt}[Stochastic oracle]
\label{ass_noise}
For every \(n\geq1\), the stochastic gradient \(G_n\) satisfies the
following conditions.
\begin{enumerate}[label=\textnormal{(\roman*)},leftmargin=*]

    \item
    \label{ass_noise:i}
    \textbf{Conditional unbiasedness.}
    \[
        \E[G_n\mid\mathscr F_{n-1}]
        =
        \nabla g(\theta_n).
    \]

    \item
    \label{ass_noise:i2}
    \textbf{Affine conditional second moment.}
    There exist constants \(\beta\geq0\) and \(\sigma\geq0\) such that
    \[
        \E[
            \|G_n\|^2
            \mid
            \mathscr F_{n-1}
        ]
        \leq
        \beta\|\nabla g(\theta_n)\|^2+\sigma^2.
    \]
\end{enumerate}
\end{assumpt}

The affine second-moment bound allows the stochastic-gradient norm to grow
with the true gradient and imposes no almost-sure uniform bound on
\(\|G_n\|\).  Together with conditional unbiasedness, it is the only oracle
condition used in the stability analysis.

Throughout the paper, we write the oracle noise as
\begin{equation}
\label{eq:oracle-noise}
    \varepsilon_n:=G_n-\nabla g(\theta_n).
\end{equation}
Under conditional unbiasedness,
\(\E[\varepsilon_n\mid\mathscr F_{n-1}]=0\).

Because \(S_n=S_{n-1}+\|G_n\|^2\) contains the current stochastic
gradient, the effective stepsize \(\alpha_0/\sqrt{S_n}\) is
\(\mathscr F_n\)-measurable but is generally not
\(\mathscr F_{n-1}\)-measurable. Consequently, it cannot be taken
outside a conditional expectation given \(\mathscr F_{n-1}\), and the
usual predictable-stepsize martingale decomposition does not apply
directly.

\subsection{Stationarity Assumptions and Their Scope}
\label{subsec:asymptotic-assumptions}

For stationarity we use a localized higher-to-second moment comparison,
requiring one fixed moment of order \(2+\rho\) only when the true gradient is
small. Related conditional higher-moment assumptions appear in martingale limit
theory \citep{machkouri2007exact,fan2019exact}.

\begin{assumpt}[Local conditional relative-moment condition]
\label{ass:local-relative-moment}
There exist \(\delta_0,\rho>0\) and \(\kappa_\rho<\infty\) such that,
with \(\mathcal{L}_n:=\{\|\nabla g(\theta_n)\|\leq\delta_0\}\), almost surely
for every \(n\),
\begin{equation}
\label{eq:local-relative-moment}
 \ind_{\mathcal{L}_n}\E[\|G_n\|^{2+\rho}\mid\mathscr F_{n-1}]
 \leq \kappa_\rho\ind_{\mathcal{L}_n}\E[\|G_n\|^2\mid\mathscr F_{n-1}].
\end{equation}
\end{assumpt}

On \(\mathcal L_n\), \cref{eq:local-relative-moment} bounds the fraction of
the conditional second moment contributed by \(\|G_n\|>R\) by
\(\kappa_\rho R^{-\rho}\). Local almost-sure boundedness implies the condition,
but the converse need not hold; in particular, locally nondegenerate additive
Gaussian noise is allowed. An ordinary uniform local \((2+\rho)\)-moment
bound is insufficient for the general adapted-oracle model; see
Proposition~\ref{prop:ordinary-local-moment-insufficient}.

\begin{assumpt}[Weak Sard condition]
\label{extra:i2}
The critical-value set \(g(\Theta_{\mathrm{stat}})\) is nowhere dense in \(\R\).
\end{assumpt}

The weak Sard condition is not used in the stability argument. It enters only
in the ODE limit-set identification; see, e.g., \citet{benaim2006dynamics}.
In particular, it allows internally chain-transitive limit sets of the
gradient-flow ODE \(\dot z=-\nabla g(z)\) to be identified with subsets of
the stationary set.

\paragraph{Examples and scope of the local relative-moment condition.}
The following propositions verify the condition for fixed-size finite-sum
mini-batches and for locally nondegenerate additive noise with a finite
\((2+\rho)\) moment. The third proposition delineates the scope of this
condition in the general adapted-oracle model.

\begin{pros}[Fixed-size mini-batch gradients]
\label{prop:ass:sto:bound}
Let $g(\theta)=\frac1N\sum_{i=1}^N g_i(\theta),$
where each \(g_i:\R^d\to\R\) is continuously differentiable. Suppose
that \cref{ass_g_poi}~\ref{ass_g_poi:i3} holds. Then every stochastic
gradient obtained by uniformly sampling a mini-batch of a fixed size,
with or without replacement, satisfies
\cref{ass:local-relative-moment} for every \(\rho>0\).
\end{pros}

\begin{pros}[Additive unbounded noise]
\label{prop:additive-noise-lrm}
Let \(G_n=\nabla g(\theta_n)+\varepsilon_n\), where
\(\E[\varepsilon_n\mid\mathscr F_{n-1}]=0\).  Suppose that on the local
region there are constants \(M_{2+\rho}<\infty\) and \(v_0>0\) such that
\(\E[\|\varepsilon_n\|^{2+\rho}\mid\mathscr F_{n-1}]\leq M_{2+\rho}\)
and \(\E[\|\varepsilon_n\|^2\mid\mathscr F_{n-1}]\geq v_0\).
Then \cref{ass:local-relative-moment} holds.  In particular, the condition
allows additive Gaussian noise with uniformly bounded covariance and a
positive local trace lower bound.
\end{pros}

\begin{pros}[An ordinary local higher-moment bound is insufficient]
\label{prop:ordinary-local-moment-insufficient}
For every \(p=2+\rho>2\), there is a one-dimensional smooth,
nonnegative, coercive objective with the non-flatness and weak Sard
properties.  It admits an adapted conditionally unbiased oracle with an
affine conditional second moment and a uniformly bounded conditional
\(p\)-moment on \(\{|\nabla g(\theta_n)|\leq\delta_0\}\), for which the
AdaGrad-Norm iterates are bounded but
\(|\nabla g(\theta_n)|\not\to0\) almost surely.
\end{pros}

The construction in \cref{app:ordinary-local-moment-counterexample} uses
increasingly rare oracle values comparable with the accumulated normalizer
scale.  Its conditional law is allowed to depend on the past through
\(S_{n-1}\), as permitted by the adapted-oracle model used here.  Thus the
counterexample does not rule out sufficiency of an ordinary local
\((2+\rho)\)-moment bound for a fresh i.i.d. state-only oracle.  Its role is
to separate the general adapted setting from stronger oracle models.  The
relative-moment condition excludes this behavior while retaining unbounded
finite-\((2+\rho)\)-moment noise.  Infinite-variance noise remains outside
the affine second-moment framework.

\section{Main Results and Proof Strategy}
\label{sec:main-results}

We now state the main result and summarize the two-stage proof architecture
used in \cref{sec:main-analysis}.

\subsection{Main Stability and Stationarity Theorem}
\label{subsec:main-theorem}

\begin{thm}[Main stability and stationarity theorem]
\label{thm:main}
Let \(\{(\theta_n,S_n)\}_{n\geq1}\) be generated by
\cref{eq:adagrad-norm-intro}, with deterministic
\(\theta_1\in\R^d\), \(S_0>0\), and \(\alpha_0>0\).
\begin{enumerate}[label=\textnormal{(\roman*)},leftmargin=*]

\item
\label{thm:main:stability}
Suppose that \cref{ass_g_poi} and
\cref{ass_noise}~\ref{ass_noise:i}--\ref{ass_noise:i2} hold.
Then the running objective supremum is integrable:
\begin{equation}
\label{eq:main-expected-supremum}
    \E\!\left[\sup_{n\geq1} g(\theta_n)\right]<\infty.
\end{equation}

\item
\label{thm:main:boundedness}
Under the assumptions of part~\ref{thm:main:stability}, the iterates are
integrably bounded:
\begin{equation}
\label{eq:main-iterate-boundedness}
    \E\!\left[\sup_{n\geq1}\|\theta_n\|\right]<\infty.
\end{equation}
In particular, $ \sup_{n\geq1}\|\theta_n\|<\infty
    \qquad\text{almost surely}.$

\item
\label{thm:main:stationarity}
If, in addition,
\cref{ass:local-relative-moment,extra:i2} hold, then
\begin{equation}
\label{eq:main-stationarity}
    \|\nabla g(\theta_n)\|\to0
    \quad\text{almost surely},
    \qquad
    \E\|\nabla g(\theta_n)\|^2\to0.
\end{equation}

\end{enumerate}
\end{thm}

\subsection{Interpretation and Proof Strategy}
\label{subsec:interpretation-proof}
\begin{figure}[t]
\centering
\begin{tikzpicture}[
    font=\footnotesize,
    box/.style={
        draw,
        rounded corners=2pt,
        align=center,
        inner xsep=4pt,
        inner ysep=3pt,
        minimum height=18mm
    },
    sone/.style={
        box,
        text width=0.19\textwidth
    },
    stwoleft/.style={
        box,
        text width=0.23\textwidth
    },
    stwomid/.style={
        box,
        text width=0.285\textwidth
    },
    stworight/.style={
        box,
        text width=0.23\textwidth
    },
    arrow/.style={
        -{Latex[length=1.6mm]},
        line width=0.45pt
    },
    stageframe/.style={
        draw,
        dashed,
        rounded corners=3pt,
        line width=0.4pt,
        inner xsep=6pt,
        inner ysep=6pt
    },
    stagelabel/.style={
        font=\bfseries,
        fill=white,
        inner xsep=3pt,
        inner ysep=0pt
    }
]

%================================================
% Stage I: Stability
%================================================

\node[sone] (stabass)
{\textbf{Assumptions}\\
smoothness, non-flatness,\\
unbiasedness, affine second moment};

\node[sone, right=4mm of stabass] (lyap)
{\textbf{Current-sample}\\
\textbf{Lyapunov correction}\\
predictable descent};

\node[sone, right=4mm of lyap] (exc)
{\textbf{Three-phase}\\
\textbf{stopped excursion}\\
controls nonsummable residuals};

\node[sone, right=4mm of exc] (stab)
{\textbf{Integrable stability}\\[0.5mm]
$\E[\sup_n g(\theta_n)]<\infty$\\
$\E[\sup_n\|\theta_n\|]<\infty$};

\draw[arrow] (stabass) -- (lyap);
\draw[arrow] (lyap) -- (exc);
\draw[arrow] (exc) -- (stab);

\node[
    stageframe,
    fit=(stabass)(lyap)(exc)(stab)
] (stage1) {};

\node[
    stagelabel,
    anchor=west
] at ($(stage1.north west)+(5mm,0)$)
{Stage I: Stability};

%================================================
% Stage II: Stationarity
%================================================

\node[
    stwomid,
    below=12mm of stage1.south
] (regimes)
{\textbf{Two accumulator regimes}\\[0.5mm]
$\sup_n S_n<\infty$: gradient budget\\
$S_n\to\infty$: finite-window noise $+$ ODE};

\node[
    stwoleft,
    left=6mm of regimes
] (statass)
{\textbf{Additional conditions}\\
local relative-moment condition\\
and weak Sard};

\node[
    stworight,
    right=6mm of regimes
] (stat)
{\textbf{Full-sequence stationarity}\\[0.5mm]
$\|\nabla g(\theta_n)\|\to0$ a.s.\\
$\E\|\nabla g(\theta_n)\|^2\to0$};

\draw[arrow] (statass) -- (regimes);
\draw[arrow] (regimes) -- (stat);

\node[
    stageframe,
    fit=(statass)(regimes)(stat)
] (stage2) {};

\node[
    stagelabel,
    anchor=west
] at ($(stage2.north west)+(5mm,0)$)
{Stage II: Stationarity};

%================================================
% Connection between stages
%================================================

\draw[arrow]
    (stage1.south)
    --
    (stage2.north);

\end{tikzpicture}

\caption{
Proof architecture for Theorem~\ref{thm:main}.
Stage~I establishes integrable trajectory stability.
Stage~II combines this stability with the additional stationarity
conditions and treats the two accumulator regimes separately to obtain
full-sequence stationarity.
}
\label{fig:proof-architecture}
\end{figure}
Theorem~\ref{thm:main} has two layers.
Parts~\ref{thm:main:stability}--\ref{thm:main:boundedness} establish
trajectory stability under the objective and oracle assumptions.  Iterate
boundedness is derived rather than assumed: the objective running supremum is
controlled first, and the linear-growth consequence of non-flatness then gives
integrable iterate boundedness.  Neither the local relative-moment condition
nor weak Sard is used at this stage.
Part~\ref{thm:main:stationarity} adds these two conditions for asymptotic
stationarity.  It gives full-sequence gradient-norm convergence almost surely
and in mean square, without asserting convergence of \(\theta_n\) to a single
stationary point.

The proof follows the two-stage architecture in
\cref{fig:proof-architecture}.  Stage~I combines a current-sample Lyapunov
correction with a three-phase stopped-excursion argument to establish
integrable trajectory stability.  Stage~II treats bounded accumulators
through the normalized gradient budget and divergent accumulators through
finite-window noise control and the gradient-flow ODE.  Integrable stability
then upgrades almost-sure stationarity to mean-square convergence.

\section{Stability and Stationarity Analysis of AdaGrad-Norm}
\label{sec:main-analysis}

We now prove \cref{thm:main} for the standard square-root recursion
\cref{eq:adagrad-norm-intro}; power-normalized extensions are deferred to
\cref{sec:scalar-extensions}.

\subsection{Integrable Trajectory Stability}
\label{sec:stability}

We first prove parts~\ref{thm:main:stability}--\ref{thm:main:boundedness};
the iterate bound then follows from \cref{lem:nonflat-linear-growth}.

\subsubsection{The Current-Sample Lyapunov Structure}
\label{subsec:current-sample-lyapunov}

For square-root AdaGrad-Norm, define
\begin{equation}
\label{eq:adagrad-normalized-quantities}
    \zeta_n
    :=
    \frac{\|\nabla g(\theta_n)\|^2}{\sqrt{S_{n-1}}},
    \qquad
    \Gamma_n
    :=
    \frac{\|G_n\|^2}{S_n},
    \qquad
    \widetilde\Gamma_n
    :=
    \frac{\Gamma_n}{\sqrt{S_n}}.
\end{equation}
Here, \(\zeta_n\) is the predictable descent quantity,
\(\Gamma_n\) is the critical square-root residual, and
\(\widetilde\Gamma_n\) is the higher-order residual. We use the adapted Lyapunov process
\begin{equation}
\label{eq:adagrad-lyapunov}
    V_n
    :=
    g(\theta_n)
    +
    \frac{\beta\alpha_0}{2}\zeta_n.
\end{equation}
Since \(S_{n-1}\geq S_0\) and
\(\|\nabla g(\theta_n)\|^2\leq2Lg(\theta_n)\),
\begin{equation}
\label{eq:adagrad-lyapunov-comparison}
    g(\theta_n)
    \leq
    V_n
    \leq
    \kappa_Vg(\theta_n),
    \qquad
    \kappa_V
    :=
    1+
    \frac{\beta\alpha_0L}{\sqrt{S_0}}.
\end{equation}

\paragraph{Adapted Lyapunov descent.}

\begin{lem}[Sufficient descent inequality]
\label{sufficient:lem}
Suppose that
\cref{ass_g_poi}~\ref{ass_g_poi:i2} and
\cref{ass_noise}~\ref{ass_noise:i}--\ref{ass_noise:i2} hold.
Then there exist constants
\(C_{\Gamma,1},C_{\Gamma,2}<\infty\) and a
martingale-difference sequence
\(\{\widehat X_n\}_{n\geq1}\) such that
\begin{equation}
\label{eq:sufficient-descent}
    V_{n+1}-V_n
    \leq
    -\frac{\alpha_0}{4}\zeta_n
    +
    C_{\Gamma,1}\Gamma_n
    +
    C_{\Gamma,2}\widetilde\Gamma_n
    +
    \alpha_0\widehat X_n,
\end{equation}
where \(\E[\widehat X_n\mid\mathscr F_{n-1}]=0\).
\end{lem}

\begin{namedproof}{Proof idea.}
Rewrite the current denominator around \(S_{n-1}^{-1/2}\), center the
resulting correction, and absorb its leading conditional mean with the
Lyapunov term \((\beta\alpha_0/2)\zeta_n\).  The remaining terms are
controlled by \(\Gamma_n\) and \(\widetilde\Gamma_n\); see
\cref{proof:sufficient:lem}.
\end{namedproof}

For later use, we associate \cref{eq:sufficient-descent} with
\[
(D_n,R_n,\widetilde R_n,M_n)
=(\zeta_n,\Gamma_n,\widetilde\Gamma_n,\alpha_0\widehat X_n).
\]

The inequality in \cref{eq:sufficient-descent} is the basic structural
estimate for the critical case.  Its negative term is predictable, but the
residual \(\Gamma_n\) is only logarithmically controlled globally.  The role
of the stopped-excursion argument is to use this descent only where it is
needed: martingale variation is controlled on a bounded entrance band, while
accumulated descent controls the subsequent return phase.

\subsubsection{A Three-Phase Stopped-Excursion Principle}
\label{sec:abstract-framework}

To convert the preceding Lyapunov inequality into a running-supremum bound,
we use the following abstract stopped-excursion principle.  We retain the
notation \(V_n\) so that its later specialization to
\cref{eq:adagrad-lyapunov} is transparent.

Let \(\{V_n\}_{n\geq1}\) be a nonnegative process such that \(V_n\) is
\(\mathscr F_{n-1}\)-measurable, and let
\(\{Y_n\}_{n\geq1}\) be another nonnegative
\(\mathscr F_{n-1}\)-measurable process satisfying
\(Y_n\leq V_n\). Suppose that
\begin{equation}
\label{eq:abstract-drift}
    V_{n+1}-V_n
    \leq
    -c_DD_n
    +
    C_RR_n
    +
    C_H\widetilde R_n
    +
    M_n,
\end{equation}
where \(c_D,C_R,C_H>0\), \(D_n\geq0\) is
\(\mathscr F_{n-1}\)-measurable,
\(R_n,\widetilde R_n\geq0\) are integrable and
\(\mathscr F_n\)-measurable, and \(M_n\) is square integrable,
\(\mathscr F_n\)-measurable, and conditionally mean zero given
\(\mathscr F_{n-1}\).

Fix a deterministic level \(G>0\) such that \(V_1\leq G\), and put
\(c_0:=0\).  With \(\inf\varnothing:=\infty\), define recursively
\begin{align}
    a_i
    &:=
    \inf\{n>c_{i-1}:V_n>G\},
    \label{eq:three-phase-entrance}\\
    b_i
    &:=
    \inf\{n\geq a_i:V_n\leq G\ \text{or}\ V_n>2G\},
    \label{eq:three-phase-band-exit}\\
    c_i
    &:=
    \inf\{n\geq b_i:V_n\leq G\}.
    \label{eq:three-phase-return}
\end{align}
If one of these sets is empty, that time and all subsequent times are
set equal to \(\infty\).
The three phases are \(I_i^1=[a_i,b_i)\), \(I_i^2=[b_i,c_i)\), and
\(I_i^3=[c_i,a_{i+1})\).
The first phase lies in the band \(G<V_n\leq2G\).  The second is
nonempty only after an upcrossing of \(2G\) and continues until the
process returns to \((-\infty,G]\).  The third lies below \(G\).
Since \(V_n\) is \(\mathscr F_{n-1}\)-measurable, the selectors
\begin{equation}
\label{eq:three-phase-selectors}
    H^1_{i,n}:=\ind_{\{a_i\leq n<b_i\}},
    \qquad
    H^2_{i,n}:=\ind_{\{b_i\leq n<c_i\}}
\end{equation}
are \(\mathscr F_{n-1}\)-measurable.  The intervals are defined on
the infinite horizon; deterministic horizons are introduced only by
restricting sums to \(n<T\).
\begin{lem}[Predictability of the phase selectors]
\label{lem:three-phase-predictability}
For every \(i,n\geq1\), the events \(\{a_i=n\}\), \(\{b_i=n\}\), and
\(\{c_i=n\}\) belong to \(\mathscr F_{n-1}\).
Consequently, \(H^1_{i,n}\) and \(H^2_{i,n}\) are
\(\mathscr F_{n-1}\)-measurable, and
\(\{H^r_{i,n}M_n\}_{n\geq1}\) is a martingale-difference sequence
for \(r\in\{1,2\}\).
\end{lem}
The measurability claim follows directly from the finite-history definition
of the phase times; the complete argument is given in
Appendix~\ref{app:three-phase-proof}.

\begin{thm}[Three-phase stopped-excursion criterion]
\label{thm:abstract-stability}
For the deterministic level \(G\) and the phase times introduced above,
suppose that \cref{eq:abstract-drift} holds.  Suppose also that there
are deterministic constants \(B_G,B_{2G}>0\), with \(B_G<2G\), and
finite nonnegative constants
\(\kappa_0,\kappa_1,\kappa_2,\kappa_3\) satisfying the following four
verification conditions.

\begin{enumerate}[label=\textnormal{(\roman*)},leftmargin=*]

    \item
    \textbf{Controlled crossings.}
    For every \(i\),
    \[
        V_{a_i}\leq B_G
        \quad\text{on }\{a_i<\infty\},
        \qquad
        V_{b_i}\leq B_{2G}
        \quad\text{on }\{b_i<\infty,\ V_{b_i}>2G\}.
    \]

    \item
    \textbf{Residual budgets above the Lyapunov threshold.}
    \begin{equation}
    \label{eq:abstract-uniform-costs}
    A_R
    :=
    \sum_{n\geq1}\E[\ind_{\{V_n>G\}}R_n]<\infty,
    \qquad
    A_H
    :=
    \sum_{n\geq1}\E[\ind_{\{V_n>G\}}\widetilde R_n]<\infty.
    \end{equation}

    \item
    \textbf{Conditional quadratic variation on the entrance band.}\par
    For every \(n\geq1\),
    \begin{equation}
    \label{eq:abstract-martingale-control}
    \begin{aligned}
        \ind_{\{G<V_n\leq2G\}}
        \E[M_n^2\mid\mathscr F_{n-1}]
        &\leq
        \ind_{\{G<V_n\leq2G\}}
        \left[
          (\kappa_0+\kappa_1V_n+\kappa_2V_n^2)
          \E[R_n\mid\mathscr F_{n-1}]
          \right.\\[-1mm]
        &\left.
          +\kappa_3\E[\widetilde R_n\mid\mathscr F_{n-1}]
        \right].
    \end{aligned}
    \end{equation}

    \item
    \textbf{Target variation above the Lyapunov threshold.}\par
    There are finite constants \(C_{Y,D},C_{Y,R},C_{Y,H}\) such that
    \begin{equation}
    \label{eq:abstract-target-variation}
    \begin{aligned}
        \ind_{\{V_n>G\}}
        \E[|Y_{n+1}-Y_n|\mid\mathscr F_{n-1}]
       & \leq
        \ind_{\{V_n>G\}}
        \left(
          C_{Y,D}D_n
          +C_{Y,R}\E[R_n\mid\mathscr F_{n-1}]
          \right.\\[-1mm]
        & \left.
          +C_{Y,H}\E[\widetilde R_n\mid\mathscr F_{n-1}]
        \right).
    \end{aligned}
    \end{equation}

\end{enumerate}

Then, by monotone convergence,
\begin{equation}
\label{eq:abstract-integrable-stability}
 \E\!\left[\sup_{n\geq1}Y_n\right]
 =\sup_{T\geq1}\E\!\left[\max_{1\leq n\leq T}Y_n\right]<\infty.
\end{equation}
In particular, \(\sup_{n\geq1}Y_n<\infty\) almost surely, and
\(\{Y_n\}_{n\geq1}\) is uniformly integrable.
\end{thm}

\begin{namedproof}{Proof sketch.}
Phase~I controls upcrossings of \(2G\) through the residual and local
martingale budgets; Phase~II uses negative drift and target variation to
control the return excursion; Phase~III stays below \(G\).  Summing the
disjoint excursions gives the uniform \(L^1\) envelope; see
Appendix~\ref{app:three-phase-proof}.
\end{namedproof}

The criterion localizes martingale variation to the bounded entrance band and
uses negative drift only on the return phase; no global quadratic-variation
bound is required.

\subsubsection{Verification for AdaGrad-Norm}
\label{subsec:adagrad-verification}

We now connect the concrete inequality in \cref{eq:sufficient-descent} to
\cref{thm:abstract-stability}.  For readability, we organize the verification
into five steps, \textnormal{(N1)}--\textnormal{(N5)}.  Step \textnormal{(N1)}
is precisely Lemma~\ref{sufficient:lem}; the remaining steps establish the
crossing, residual, martingale-variation, and target-variation estimates
required by the three-phase criterion.

\begin{table}[t]
\centering
\caption{Verification map for square-root AdaGrad-Norm.  The five steps
connect the concrete Lyapunov inequality to the hypotheses of the
three-phase stopped-excursion criterion.}
\label{tab:adagrad-verification-map}
\renewcommand{\arraystretch}{1.13}
\setlength{\tabcolsep}{2.5pt}
\small
\begin{tabularx}{\textwidth}{@{}
>{\centering\arraybackslash}p{0.10\textwidth}
>{\raggedright\arraybackslash}p{0.31\textwidth}
>{\raggedright\arraybackslash}X@{}}
\toprule
Step & AdaGrad-Norm result & Role in the stopped-excursion argument \\
\midrule
N1 & Lemma~\ref{sufficient:lem}
   & Adapted Lyapunov descent after the current-sample correction \\
N2 & Lemmas~\ref{lem:adj:ghat} and~\ref{pro_0}
   & Threshold separation and controlled crossings \\
N3 & Lemmas~\ref{inequ:gamma:sn} and~\ref{lem_su}
   & Higher-order and high-gradient residual budgets \\
N4 & Lemma~\ref{lem:xhat}
   & Conditional quadratic variation on the entrance band \\
N5 & Lemma~\ref{lem:high-level-objective-variation} with~\ref{pro_0}
   & Objective variation on the return phase \\
\bottomrule
\end{tabularx}
\end{table}

The steps are logically coupled rather than independent.  In particular,
\textnormal{(N3)} uses the negative drift from \textnormal{(N1)}, while the
threshold separation in \textnormal{(N2)} transfers the high-gradient bounds
from \textnormal{(N3)} and \textnormal{(N5)} to the region \(\{V_n>G\}\).

\paragraph{N2: Lyapunov-Threshold Separation and Crossing Control.}

We begin by establishing a deterministic one-step growth estimate.

\begin{lem}[One-step growth of the Lyapunov process]
\label{lem:adj:ghat}
There exist constants \(a_h,b_h,C_0>0\) such that
\begin{equation}
\label{eq:lyapunov-one-step-growth}
    V_{n+1}-V_n
    \leq
    h(V_n),
    \qquad
    h(x):=a_h\sqrt{x}+b_h,
\end{equation}
and \(h(x)<x/2\) for every \(x\geq C_0\).
\end{lem}

See Appendix~\ref{app:adagrad-verification-details} for the proof.

\begin{lem}[Separation and controlled crossings]
\label{pro_0}
Suppose that
\cref{ass_g_poi}~\ref{ass_g_poi:i3} holds. Then there exist
\(\delta>0\) and \(G_0<\infty\) such that, for every
\(G\geq G_0\),

\begin{enumerate}[label=\textnormal{(\roman*)},leftmargin=*,itemsep=1pt]

    \item
    \(V_n>G\) implies \(\|\nabla g(\theta_n)\|>\delta\);

    \item if \(V_{n-1}\leq G<V_n\), then
    \(V_n\leq G+h(G)<3G/2\);

    \item if \(V_{n-1}\leq2G<V_n\), then
    \(V_n\leq2G+h(2G)\).

\end{enumerate}
\end{lem}

Lemmas~\ref{lem:adj:ghat} and~\ref{pro_0} verify step
\textnormal{(N2)}; proofs are in Appendix~\ref{app:adagrad-verification-details}.

\paragraph{N3: Self-Normalized Residual Budgets.}

\begin{lem}[Summability of the higher-order residual]
\label{inequ:gamma:sn}
The higher-order residual satisfies
\begin{equation}
\label{eq:adagrad-higher-order-residual}
    \sum_{n=1}^{\infty}\widetilde\Gamma_n
    \leq
    \frac{2}{\sqrt{S_0}}
    \quad\text{a.s.},
    \qquad
    \sum_{n=1}^{\infty}\E[\widetilde\Gamma_n]
    \leq
    \frac{2}{\sqrt{S_0}}.
\end{equation}
\end{lem}
\begin{proof}
This is the series--integral estimate
\(\sum_n(S_n-S_{n-1})S_n^{-3/2}\leq
\int_{S_0}^{\infty}x^{-3/2}\,\mathrm dx\); see
Appendix~\ref{app:adagrad-verification-details}.
\end{proof}

\begin{lem}[Normalized gradient budget away from stationarity]
\label{lem_su}
Suppose that
\cref{ass_g_poi}~\ref{ass_g_poi:i2} and
\cref{ass_noise}~\ref{ass_noise:i}--\ref{ass_noise:i2} hold.
Then, for every \(\delta>0\),
\begin{equation}
\label{eq:adagrad-high-gradient-energy}
  \sum_{n=1}^{\infty}
    \E\left[
        \ind_{\{
            \|\nabla g(\theta_n)\|>\delta
        \}}
        \Gamma_n
    \right]
    \leq
    \left(
        \beta+\frac{\sigma^2}{\delta^2}
    \right)M_\delta
\end{equation}
where \(M_\delta<\infty\) depends only on problem parameters,
algorithmic parameters, initialization, and \(\delta\).
\end{lem}

The proof is given in \cref{sec:proof:lemsu}.

Lemmas~\ref{inequ:gamma:sn} and~\ref{lem_su} complete
verification step \textnormal{(N3)}.  For every sufficiently large \(G\) supplied by
\textnormal{(N2)}, its separation property transfers the high-gradient
estimate in \cref{eq:adagrad-high-gradient-energy} to the region
\(\{V_n>G\}\); \cref{eq:adagrad-higher-order-residual} supplies the global
\(\widetilde\Gamma_n\)-budget.  These are the residual estimates required by
the stopped-excursion criterion.

\paragraph{N4: Conditional Quadratic Variation on the Entrance Band.}

Let \(\widehat X_n\) be the martingale-difference term in
\cref{sufficient:lem}; its explicit decomposition is given in
Appendix~\ref{proof:sufficient:lem}. Put
\(a_\Gamma:=\sigma^2/(2\sqrt{S_0})\).

\begin{lem}[Conditional martingale-variation bound]
\label{lem:xhat}
The process \(\widehat X_n\) is square integrable at every fixed index and
satisfies
\begin{equation}
\label{eq:xhat-conditional-variance}
 \E[(\alpha_0\widehat X_n)^2\mid\mathscr F_{n-1}]
 \leq(3\alpha_0^2a_\Gamma^2+6L\alpha_0^2V_n+3V_n^2)
 \E[\Gamma_n\mid\mathscr F_{n-1}].
\end{equation}
Hence, on \(G<V_n\leq2G\), the coefficient is at most
\(3\alpha_0^2a_\Gamma^2+12L\alpha_0^2G+12G^2\).
\end{lem}

See Appendix~\ref{app:adagrad-verification-details} for the proof.

\paragraph{N5: Objective Variation Above the Lyapunov Threshold.}

\begin{lem}[Objective variation above the Lyapunov threshold]
\label{lem:high-level-objective-variation}
For every \(\delta>0\), there is a finite constant
\(C_{\mathrm{var}}(\delta)\), depending only on \(\delta\) and the model
and algorithmic parameters, such that
\begin{equation}
\label{eq:adagrad-high-level-variation}
\begin{aligned}
    &\ind_{\{\|\nabla g(\theta_n)\|>\delta\}}
    \E[
      |g(\theta_{n+1})-g(\theta_n)|
      \mid\mathscr F_{n-1}
    ]
    \leq
    C_{\mathrm{var}}(\delta)
    \ind_{\{\|\nabla g(\theta_n)\|>\delta\}}
    \zeta_n.
\end{aligned}
\end{equation}
\end{lem}

See Appendix~\ref{app:adagrad-verification-details} for the proof.

By \cref{pro_0}, for every sufficiently large \(G\) the event
\(\{V_n>G\}\) is contained in
\(\{\|\nabla g(\theta_n)\|>\delta\}\).  Hence
\cref{lem:high-level-objective-variation} completes verification step \textnormal{(N5)} with
\(C_{5,D}=C_{\mathrm{var}}(\delta)\) and \(C_{5,R}=C_{5,H}=0\).

\subsubsection{Proof of Integrable Trajectory Stability}

\begin{namedproof}{Proof of Theorem~\ref{thm:main}, parts~\ref{thm:main:stability}--\ref{thm:main:boundedness}.}
Set
\[
(D_n,R_n,\widetilde R_n,M_n)
=(\zeta_n,\Gamma_n,\widetilde\Gamma_n,\alpha_0\widehat X_n),
\qquad Y_n:=g(\theta_n).
\]
Lemma~\ref{sufficient:lem} gives \cref{eq:abstract-drift}. Choose \(G\) large
enough that \(V_1\leq G\), \(h(G)<G/2\), and the separation in
\cref{pro_0} holds. Lemmas~\ref{lem:adj:ghat} and~\ref{pro_0} give the crossing
bounds with \(B_G=G+h(G)<3G/2\) and \(B_{2G}=2G+h(2G)\).
The same separation transfers the high-gradient \(\Gamma_n\)-budget in
\cref{lem_su} to \(\{V_n>G\}\), while \cref{inequ:gamma:sn} supplies the
global \(\widetilde\Gamma_n\)-budget. Lemmas~\ref{lem:xhat} and
\ref{lem:high-level-objective-variation} verify the entrance-band quadratic
variation and return-phase target variation. Thus
\cref{thm:abstract-stability} yields
\cref{eq:main-expected-supremum}.

Finally, \cref{lem:nonflat-linear-growth} gives, pathwise,
\(
 \sup_n\|\theta_n\|
 \leq R_c+c^{-1}\sup_n g(\theta_n)
\)
for any \(c\in(0,\ell_\infty)\). Taking expectations proves
\cref{eq:main-iterate-boundedness}, and almost-sure boundedness follows.
\end{namedproof}

This completes the stability layer of \cref{thm:main}.

\subsection{Asymptotic Stationarity}
\label{sec:asympt:result}

It remains to prove part~\ref{thm:main:stationarity}.  Since \(S_n\) is
nondecreasing, every trajectory belongs to one of the two regimes
\[
\mathcal A:=\{\lim_nS_n<\infty\},
\qquad
\mathcal A^c:=\{S_n\to\infty\}.
\]
The first is handled by the normalized gradient budget; the second by the ODE
argument below.

\subsubsection{ODE Limit and Finite-Window Noise Control}
\label{subsec:ode-identification}

We use the asymptotic-pseudotrajectory framework of
\citet{benaim2006dynamics}. Consider
\begin{equation}
\label{SA}
    z_{n+1}=z_n+\gamma_n\bigl(F(z_n)+U_n\bigr),
\end{equation}
where \(\gamma_n>0\) may be random and adapted, \(F:\R^d\to\R^d\) is locally
Lipschitz, and \(U_n\) is a stochastic perturbation. With
\(t_0:=0\), \(t_n:=\sum_{i=1}^n\gamma_i\), and
\(m(t):=\max\{j\geq0:t_j\leq t\}\), the ODE principle stated in
Appendix~\ref{app:ode-identification} requires bounded iterates,
\(\gamma_n\to0\), \(\sum_n\gamma_n=\infty\), vanishing perturbations on every
finite effective-time window, and a strict Lyapunov function whose equilibrium
values have empty interior.

For AdaGrad-Norm, \(F=-\nabla g\) is globally Lipschitz and
\begin{equation}
\label{eq:adagrad-direct-sa}
    \gamma_n:=\frac{\alpha_0}{\sqrt{S_n}},
    \qquad
    \theta_{n+1}=\theta_n+\gamma_n\bigl(-\nabla g(\theta_n)-\varepsilon_n\bigr).
\end{equation}
The current-dependent product \(\gamma_n\varepsilon_n\) is not generally a
martingale difference. Lemma~\ref{inequ:gammaU} therefore controls the actual
stepsize and proves that, for every \(T>0\),
\begin{equation}
\label{eq:adagrad-finite-window-noise}
    \lim_{n\to\infty}
    \sup_{n\leq k\leq m(t_n+T)}
    \left\|
        \sum_{i=n}^{k}
        \frac{\alpha_0}{\sqrt{S_i}}\varepsilon_i
    \right\|=0
    \qquad\text{almost surely}.
\end{equation}

\subsubsection{Almost-Sure Stationarity}
\label{subsec:almost:sure}

Under the assumptions of \cref{thm:main}~\ref{thm:main:stationarity}, on
\(\mathcal A\) the normalized gradient-budget estimate and boundedness of
\(S_n\) force every set
\(\{n:\|\nabla g(\theta_n)\|>\delta\}\) to be finite. On
\(\mathcal A^c\), iterate stability, divergence of effective time, and
\cref{inequ:gammaU} verify the ODE-identification conditions. The ODE
limit-set principle, together with weak Sard, then implies that every
accumulation point is stationary; boundedness and continuity of the gradient
yield \(\|\nabla g(\theta_n)\|\to0\) almost surely. The complete
argument is given in Appendix~\ref{appendix:asymptotic}.

\subsubsection{Mean-Square Stationarity}
\label{sec:mean:convergence}

By \cref{eq:smooth-gradient-objective-bound} and
\cref{thm:main}~\ref{thm:main:stability},
\(\|\nabla g(\theta_n)\|^2\) is dominated by the integrable random variable
\(2L\sup_k g(\theta_k)\). The almost-sure conclusion above and dominated
convergence therefore give
\(\E\|\nabla g(\theta_n)\|^2\to0\), completing
\cref{thm:main}~\ref{thm:main:stationarity}.

This completes the proof of \cref{thm:main}.  We next contrast the critical
square-root normalization with the summable power-normalized regime.

\section{Power-Normalized AdaGrad Beyond the Square-Root Boundary}
\label{sec:scalar-extensions}

We now place the critical exponent within the power-normalized family and
record the simpler summable regime above it.

\subsection{The Square-Root Boundary in the Power-Normalized Family}
\label{sec:critical-boundary}

Consider
\begin{equation}
\label{eq:power-normalized-family}
 S_n=S_{n-1}+\|G_n\|^2,
 \qquad \theta_{n+1}=\theta_n-\alpha_0S_n^{-q}G_n,
 \quad q>0.
\end{equation}
Standard AdaGrad-Norm corresponds to \(q=1/2\).  Since
\(x\mapsto x^{-2q}\) is decreasing,
\begin{equation}
\label{eq:power-energy-comparison}
    \sum_{n=1}^{T}\frac{\|G_n\|^2}{S_n^{2q}}
    \leq
    \int_{S_0}^{S_T}x^{-2q}\,dx
    =
    \begin{cases}
    \dfrac{S_0^{1-2q}-S_T^{1-2q}}{2q-1},&q>1/2,\\[2mm]
    \log(S_T/S_0),&q=1/2,\\[1mm]
    \dfrac{S_T^{1-2q}-S_0^{1-2q}}{1-2q},&q<1/2.
    \end{cases}
\end{equation}
Thus the quadratic normalization budget is summable for \(q>1/2\),
logarithmic at \(q=1/2\), and may grow polynomially for \(q<1/2\).
At the critical exponent the one-step movement still satisfies
\(\|\theta_{n+1}-\theta_n\|\leq\alpha_0\), whereas the same elementary
bound is unavailable below it.  Hence \(q=1/2\) is the boundary of the proof
mechanism developed here, not a universal algorithmic phase transition.

\subsection{The Summable Regime under Affine Second Moments}
\label{subsec:power-affine}

For \(q>1/2\), \cref{eq:power-energy-comparison} makes the quadratic error
summable, although the denominator still depends on the current oracle value.

\begin{thm}[Power-normalized AdaGrad under affine second moments]
\label{thm:power-affine}
Let \(g:\R^d\to[0,\infty)\) be continuously differentiable and \(L\)-smooth,
and let \cref{eq:power-normalized-family} hold with deterministic
\(\theta_1\in\R^d\), \(S_0>0\), and \(q>1/2\).  Assume, almost surely for
every \(n\),
\begin{align}
\E[G_n\mid\mathscr F_{n-1}]&=\nabla g(\theta_n),
\label{eq:power-unbiased}\\
\E[\|G_n\|^2\mid\mathscr F_{n-1}]&\leq
\beta\|\nabla g(\theta_n)\|^2+\sigma^2,
\label{eq:power-affine-moment}
\end{align}
for finite \(\beta,\sigma\geq0\).  Then \(g(\theta_n)\) converges to a finite
random variable and, almost surely,
\begin{equation}
\label{eq:power-as-stability}
 \sup_{n\geq1}g(\theta_n)<\infty,
 \qquad
 \sum_{n=1}^{\infty}
 \frac{\|\nabla g(\theta_n)\|^2}{S_{n-1}^q}<\infty.
\end{equation}
If \(q\leq1\), then
\begin{equation}
\label{eq:power-as-stationarity}
 \|\nabla g(\theta_n)\|\to0\qquad\text{almost surely}.
\end{equation}
This stationarity conclusion requires neither coercivity nor a condition on
critical values, and no conditional moment above order two.
\end{thm}

The restriction \(q\leq1\) in the stationarity claim reflects the need for
divergent effective optimization time; no stationarity claim is made here for
\(q>1\).

\subsection{An Integrable Upgrade under Uniform Second Moments}
\label{subsec:power-uniform-l2}

A uniform conditional second-moment bound upgrades the preceding pathwise
stability to an integrable running-supremum bound.

\begin{thm}[Integrable stability under a uniform second-moment bound]
\label{thm:power-uniform-l2}
In the setting of \cref{thm:power-affine}, suppose further that, for some
deterministic \(K<\infty\),
\begin{equation}
\label{eq:power-uniform-second-moment}
 \E[\|G_n\|^2\mid\mathscr F_{n-1}]\leq K
 \qquad\text{almost surely for every }n.
\end{equation}
Then
\begin{equation}
\label{eq:power-integrable-stability}
 \E\!\left[\sup_{n\geq1}g(\theta_n)\right]<\infty.
\end{equation}
If \(q\leq1\), then also
\begin{equation}
\label{eq:power-ms-stationarity}
 \E\|\nabla g(\theta_n)\|^2\to0.
\end{equation}
Thus a uniform conditional \(L^2\) bound, rather than an almost-sure bound on
stochastic gradients, suffices for integrable stability and mean-square
stationarity in the summable regime.
\end{thm}

Proofs of \cref{thm:power-affine,thm:power-uniform-l2} are given in
Appendix~\ref{app:power-normalization-proofs}.  Fixed additive offsets in the
normalizer, or constant positive multipliers of the accumulator increment, are
exact reparameterizations obtained by changing the initial accumulator and base
stepsize; time-varying predictable weights are not.

\section{Numerical Experiments}
\label{sec:numerical}

The experiments are diagnostic rather than empirical tests of almost-sure
statements.  We examine three features suggested by the theory: unbounded
finite-moment noise at the square-root normalization, the dependence on the
normalization exponent \(q\), and the practical effect of using the current
rather than the delayed denominator.

\subsection{Experimental Protocol}
\label{subsec:numerical-protocol}

Unless stated otherwise, \(S_0=1\) and \(\alpha_0=0.5\).
We compare the current-denominator update
\(\theta_{n+1}=\theta_n-\alpha_0S_n^{-q}G_n\)
with its delayed counterpart
\(\theta_{n+1}=\theta_n-\alpha_0S_{n-1}^{-q}G_n\), where
\(S_n=S_{n-1}+\|G_n\|^2\).
Within each seed, all values of \(q\) and both denominator conventions use
the same stochastic samples.  Curves report seed medians with pointwise
\(10\%\)--\(90\%\) quantile bands; paired current--delayed comparisons use
percentile bootstrap intervals based on \(10{,}000\) resamples.

We monitor the checkpointed tail gradient
\(\mathcal G_T^{\rm tail}
:=\max\{\|\nabla g(\theta_n)\|:
n\in\mathcal C_T,\ n\geq T/2\}\),
where \(\mathcal C_T\) denotes the set of recorded checkpoints up to
horizon \(T\), and the positive objective overshoot
\(\mathcal O_T:=
[\max_{1\leq n\leq T}g(\theta_n)-g(\theta_1)]_+\).
The cumulative effective times are
\(\mathcal T_T^{\rm cur}(q)=\alpha_0\sum_{n=1}^T S_n^{-q}\) and
\(\mathcal T_T^{\rm del}(q)=\alpha_0\sum_{n=1}^T S_{n-1}^{-q}\).
The synthetic experiments use 50 seeds and \(2^{20}\) updates, while the
California Housing experiment uses 30 seeds and a budget of
\(2^{22}\) component-gradient evaluations.  Additional implementation
details are provided in the accompanying reproducibility package.

\subsection{Unbounded Finite-Moment Noise at the Square-Root Boundary}
\label{subsec:numerical-assumption-stress}

We first consider the smooth coercive nonconvex objective in dimension
\(d=20\),
\begin{equation}
\label{eq:numerical-stress-objective}
g(\theta)
=
0.025\|\theta\|^2
+
\sum_{j=1}^{20}\bigl(1-\cos(\theta_j)\bigr),
\qquad
\theta_{1,j}=(-1)^{j+1}5.5,
\end{equation}
with the conditionally unbiased oracle
\begin{equation}
\label{eq:numerical-stress-oracle}
G_n
=
\bigl(1+0.5Z_{0,n}\bigr)\nabla g(\theta_n)+Z_{1,n}u_n.
\end{equation}
Here
\(u_n=d^{-1/2}(\epsilon_{n,1},\ldots,\epsilon_{n,d})\), where the
\(\epsilon_{n,j}\) are independent Rademacher signs, so that
\(\|u_n\|=1\) deterministically.
The independent symmetric variables \(Z_{0,n}\) and \(Z_{1,n}\) follow,
in separate experiments, Rademacher, standard Gaussian,
variance-standardized Student-\(t_3\), or symmetric Pareto-\(1.5\)
distributions.  The first three satisfy the finite-moment assumptions
relevant here, whereas Pareto-\(1.5\) has infinite variance and serves only
as a stress control.

To track the perturbation appearing in the ODE analysis, let \(I_m\) denote
successive completed blocks of one unit of cumulative effective time and set
\begin{equation}
\label{eq:numerical-finite-window}
\mathcal W_m
:=
\max_{k\in I_m}
\left\|
\sum_{n=\min I_m}^{k}
\gamma_n\bigl(G_n-\nabla g(\theta_n)\bigr)
\right\|,
\qquad
\gamma_n=\frac{\alpha_0}{\sqrt{S_n}}.
\end{equation}

\begin{figure}[t]
\centering
\includegraphics[width=\textwidth]{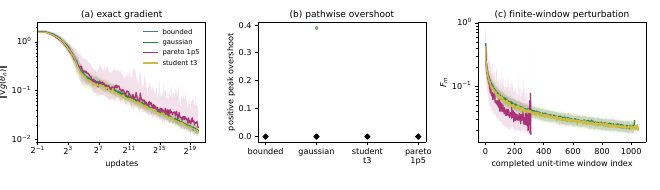}
\caption{
Assumption stress under current normalization with \(q=1/2\).
Panels report (a) the exact gradient norm, (b) objective overshoot
\(\mathcal O_T\), and (c) the finite-window perturbation
\(\mathcal W_m\).  Lines are seed medians with pointwise
\(10\%\)--\(90\%\) quantile bands; Pareto-\(1.5\) is the
infinite-variance control.
}
\label{fig:assumption-stress}
\end{figure}

\Cref{fig:assumption-stress} shows similar long-horizon behavior across
the three finite-moment noise models despite their different tail behavior,
whereas the infinite-variance control is substantially more variable.
These observations are diagnostic and do not extend the theorem beyond its
stated moment assumptions.

\subsection{Power Normalization and Effective Optimization Time}
\label{subsec:numerical-power}

We next use the Student-\(t_3\) oracle and compare
\(q\in\{0.4,0.5,0.55,0.6,0.75,1,1.1\}\) under both denominator
conventions.  The value \(q=0.4\) lies below the analyzed range and is
included only as a control.  For current normalization, the realized
normalization budget follows the three-regime comparison in
\cref{eq:power-energy-comparison}: it is uniformly bounded for
\(q>1/2\), logarithmic at \(q=1/2\), and need not remain uniformly bounded
below the square-root exponent.

\begin{figure}[t]
\centering
\includegraphics[width=\textwidth]{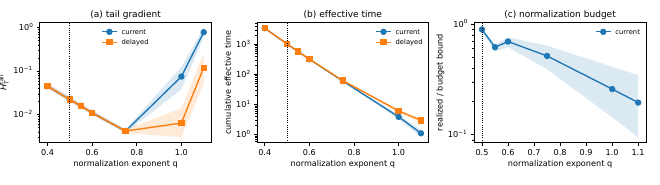}
\caption{
Power-regime diagnostics after \(2^{20}\) updates under
Student-\(t_3\) noise.  Panels report
(a) \(\mathcal G_T^{\rm tail}\),
(b) cumulative effective time, and
(c) normalization-budget usage for the current denominator.
Markers are seed medians with pointwise \(10\%\)--\(90\%\) quantile bands.
The dotted line marks \(q=1/2\), and \(q=0.4\) is an out-of-range control.
}
\label{fig:power-regimes}
\end{figure}

\Cref{fig:power-regimes} compares the tail gradient, cumulative effective
time, and normalization-budget usage.  The \(q>1\) runs illustrate that a
larger exponent can leave little effective optimization time at a fixed
computational horizon; the figure should therefore not be interpreted as a
universal ranking of the exponents or as evidence of an algorithmic phase
transition.

\subsection{Current vs.\ Delayed Normalization on a Finite-Sum Problem}
\label{subsec:numerical-california}

Finally, we consider ridge-regularized Welsch regression
\citep{barron2019general} on the California Housing data set
\citep{pace1997sparse},
\begin{equation}
\label{eq:numerical-welsch}
g_{\mathrm W}(\theta)
=
\frac1N\sum_{i=1}^{N}
c^2\left[
1-\exp\!\left(
-\frac{(x_i^\top\theta-y_i)^2}{2c^2}
\right)
\right]
+\frac{\lambda}{2}\|\theta\|^2,
\end{equation}
with \(c=1.345\) and \(\lambda=10^{-3}\).
We use a fixed training/test split, standardize the data using training
statistics, append an intercept, and initialize at \(\theta_1=0\).
Mini-batches of size 32 are sampled with replacement for
\(q\in\{0.5,0.6,0.75,1,1.1\}\).

\FloatBarrier
\begin{table}[t]
\centering
\caption{
California Housing after \(2^{22}\) component-gradient evaluations.
The two tail columns report seed medians of
\(\mathcal G_T^{\rm tail}\).
Each \(\Delta\) is the paired median difference
(current minus delayed), with a \(95\%\) percentile-bootstrap interval.
MAE is evaluated on the standardized test target.
}
\label{tab:california-current-delayed}
\renewcommand{\arraystretch}{1.12}
\scriptsize
\setlength{\tabcolsep}{3pt}
\begin{tabular}{@{}ccccc@{}}
\toprule
\(q\)
& Tail current
& Tail delayed
& \makecell{Tail \(\Delta\)\\95\% CI}
& \makecell{MAE \(\Delta\)\\95\% CI}
\\
\midrule

0.5
& 0.01583
& 0.01583
& \makecell{
  \(-1.20\mathrm{e}{-7}\)\\
  \([-2.52\mathrm{e}{-7},-2.60\mathrm{e}{-8}]\)}
& \makecell{
  \(-4.16\mathrm{e}{-9}\)\\
  \([-2.06\mathrm{e}{-8},+2.08\mathrm{e}{-9}]\)}
\\

0.6
& 0.01009
& 0.01009
& \makecell{
  \(-1.64\mathrm{e}{-7}\)\\
  \([-3.13\mathrm{e}{-7},+2.00\mathrm{e}{-8}]\)}
& \makecell{
  \(-1.83\mathrm{e}{-7}\)\\
  \([-2.57\mathrm{e}{-7},-1.56\mathrm{e}{-7}]\)}
\\

0.75
& 0.00563
& 0.00563
& \makecell{
  \(+4.42\mathrm{e}{-7}\)\\
  \([+6.43\mathrm{e}{-8},+8.93\mathrm{e}{-7}]\)}
& \makecell{
  \(-2.15\mathrm{e}{-5}\)\\
  \([-2.65\mathrm{e}{-5},-1.68\mathrm{e}{-5}]\)}
\\

1.0
& 0.00692
& 0.00689
& \makecell{
  \(+5.92\mathrm{e}{-5}\)\\
  \([+5.28\mathrm{e}{-5},+8.34\mathrm{e}{-5}]\)}
& \makecell{
  \(-1.98\mathrm{e}{-4}\)\\
  \([-2.84\mathrm{e}{-4},-1.77\mathrm{e}{-4}]\)}
\\

1.1
& 0.00944
& 0.00929
& \makecell{
  \(+1.29\mathrm{e}{-4}\)\\
  \([+8.62\mathrm{e}{-5},+1.78\mathrm{e}{-4}]\)}
& \makecell{
  \(-2.49\mathrm{e}{-4}\)\\
  \([-3.69\mathrm{e}{-4},-1.78\mathrm{e}{-4}]\)}
\\

\bottomrule
\end{tabular}
\end{table}

\Cref{tab:california-current-delayed} shows that the paired differences are
negligible for the smaller exponents and become more visible as \(q\)
increases.  Test MAE is included only as a predictive diagnostic; the
purpose of this experiment is to isolate denominator timing rather than
to rank the normalization exponents.
\section{Discussion and Conclusion}\label{sec:conclusion}
At \(q=1/2\), the current denominator and the logarithmic quadratic budget
prevent a direct application of standard stochastic-approximation arguments.
A current-sample Lyapunov correction together with the three-phase
stopped-excursion analysis yields integrable running-supremum bounds for both
the objective and the iterates.  These bounds support a two-regime ODE
argument that establishes almost-sure and mean-square convergence of the full
gradient-norm sequence.  The local conditional relative-moment condition is
needed only for the finite-window perturbation estimate and still allows
unbounded-support noise.  Our results do not cover objectives whose gradients
vanish at infinity, infinite-variance noise, finite-time rates, or convergence
to a single stationary point.  For \(q>1/2\), the relevant normalization
budget is summable, whereas for \(q<1/2\) even a uniform one-step movement
bound is unavailable from the same argument.  Weakening non-flatness at the
square-root boundary and developing useful stability conditions in the
subcritical regime remain open problems.
% AUTHOR CHECK (funding): add the full funding source and grant number that
% supported Ruinan Jin's contribution to this work. Funding support should be
% disclosed as funding rather than described as manuscript contribution.
% AUTHOR CHECK (COI): confirm that there are no additional financial competing
% interests requiring disclosure under the current JMLR policy.
\acks{The work of Xiaoyu Wang is supported in part by the Fundamental Research
Funds for the Central Universities E5EQ0101X2. The authors declare no financial
competing interests.}

\appendix

\section{Proof of the Three-Phase Stability Principle}
\label{app:three-phase-proof}

\begin{namedproof}{Proof of \cref{lem:three-phase-predictability}.}
For fixed \(i,n\), each of
\(\{a_i\leq n\}\), \(\{b_i\leq n\}\), and \(\{c_i\leq n\}\) is
determined by the finite vector \((V_1,\ldots,V_n)\).  This follows
inductively from
\cref{eq:three-phase-entrance,eq:three-phase-band-exit,eq:three-phase-return}.
Since \(V_k\) is \(\mathscr F_{k-1}\)-measurable and the filtration is
increasing, all three events belong to \(\mathscr F_{n-1}\).  Taking
the difference between the events with indices \(n\) and \(n-1\)
proves the asserted measurability of the level sets.  Moreover,
\[
 \{a_i\leq n<b_i\},\quad \{b_i\leq n<c_i\}\in\mathscr F_{n-1}.
\]
Multiplication of \(M_n\) by either selector preserves conditional
mean zero.
\end{namedproof}

\begin{namedproof}{Proof of \cref{thm:abstract-stability}.}
For \(T=1\), the conclusion follows from
\(Y_1\leq V_1\leq G\).  Fix \(T\geq2\), and let
\[
    E_i:=\{b_i<\infty,\ V_{b_i}>2G\},
    \qquad
    J_T:=\sum_{i\geq1}\ind_{E_i\cap\{b_i<T\}}.
\]
Because \(a_{i+1}>c_i\), only finitely many phase intervals intersect
\(\{1,\ldots,T\}\); all sums over \(i\) below are finite at the fixed
horizon.
Set \(\eta:=2G-B_G>0\), and, for each \(i\), define the truncated
first-phase costs
\[
    A_{i,T}
    :=
    \sum_{n<T}H^1_{i,n}(C_RR_n+C_H\widetilde R_n),
    \qquad
    Q_{i,T}
    :=
    \sum_{n<T}H^1_{i,n}M_n.
\]
On \(E_i\cap\{b_i<T\}\), summing
\cref{eq:abstract-drift} from \(a_i\) to \(b_i-1\), dropping the
nonpositive drift, and using \(V_{a_i}\leq B_G\) gives
\[
    \eta
    \leq
    A_{i,T}+|Q_{i,T}|.
\]
Consequently,
\[
    \ind_{E_i\cap\{b_i<T\}}
    \leq
    \frac{2A_{i,T}}{\eta}
    +
    \frac{4Q_{i,T}^2}{\eta^2}.
\]
The selectors in \cref{eq:three-phase-selectors} are predictable, so
martingale orthogonality and
\cref{eq:abstract-martingale-control} imply, with
\(K_G:=\kappa_0+2G\kappa_1+4G^2\kappa_2\),
\begin{align*}
    \sum_{i\geq1}\E[Q_{i,T}^2]
    &=
    \E\left[
      \sum_{i\geq1}\sum_{n<T}H^1_{i,n}
      \E[M_n^2\mid\mathscr F_{n-1}]
    \right]
    \leq K_GA_R+\kappa_3A_H.
\end{align*}
The first-phase intervals are disjoint and lie in \(\{V_n>G\}\).
It follows that
\begin{equation}
\label{eq:abstract-upcrossing-count}
    \sup_{T\geq2}\E[J_T]
    \leq
    \frac{2(C_RA_R+C_HA_H)}{\eta}
    +
    \frac{4(K_GA_R+\kappa_3A_H)}{\eta^2}
    =:C_J<\infty.
\end{equation}

The second phase is nonempty precisely after an upcrossing.  For each
such phase with \(b_i<T\), put \(e_{i,T}:=c_i\wedge T\).  Summing
\cref{eq:abstract-drift} from \(b_i\) to \(e_{i,T}-1\) gives
\[
\begin{aligned}
 c_D\sum_{n=b_i}^{e_{i,T}-1}D_n
 \leq{}&
 V_{b_i}-V_{e_{i,T}}
 +C_R\sum_{n=b_i}^{e_{i,T}-1}R_n
 +C_H\sum_{n=b_i}^{e_{i,T}-1}\widetilde R_n
 +\sum_{n=b_i}^{e_{i,T}-1}M_n .
\end{aligned}
\]
Taking expectations, summing over \(i\), using
\(V_{e_{i,T}}\geq0\), and applying
\cref{lem:three-phase-predictability} gives
\begin{equation}
\label{eq:abstract-return-energy}
\begin{aligned}
    c_D\E\left[
      \sum_{i\geq1}\sum_{n<T}H^2_{i,n}D_n
    \right]
    \leq{}&
    B_{2G}\E[J_T]
    +C_RA_R+C_HA_H.
\end{aligned}
\end{equation}
Here each initial value is at most \(B_{2G}\), and the stopped
martingale sum has expectation zero.  Thus the left-hand side of
\cref{eq:abstract-return-energy} is bounded uniformly in \(T\).

By \cref{eq:abstract-target-variation,eq:abstract-uniform-costs,eq:abstract-return-energy},
\begin{equation}
\label{eq:abstract-return-variation}
    \sup_{T\geq2}
    \E\left[
      \sum_{i\geq1}\sum_{n<T}H^2_{i,n}|Y_{n+1}-Y_n|
    \right]
    <\infty.
\end{equation}
Outside the second phases, either \(V_n\leq G\) or
\(G<V_n\leq2G\).  At the beginning of every nonempty second phase,
\(Y_{b_i}\leq V_{b_i}\leq B_{2G}\).  Hence, pathwise,
\begin{equation}
\label{eq:abstract-pathwise-envelope}
    \max_{1\leq k\leq T}Y_k
    \leq
    \max\{2G,B_{2G}\}
    +
    \sum_{i\geq1}\sum_{n<T}H^2_{i,n}|Y_{n+1}-Y_n|.
\end{equation}
The bounds in \cref{eq:abstract-return-variation,eq:abstract-pathwise-envelope}
prove the finite-horizon estimate in \cref{eq:abstract-integrable-stability}.
Monotone convergence yields the infinite-horizon bound; almost-sure
boundedness and uniform integrability follow by domination with
\(\sup_nY_n\).
\end{namedproof}

\section{Detailed Stability Verification for AdaGrad-Norm}
\label{app:adagrad-verification-details}

We give the deferred verification calculations used in
\cref{subsec:adagrad-verification}.

\begin{namedproof}{Proof of \cref{lem:xhat}.}
Since conditional centering cannot increase a conditional second
moment and \(0\leq\Lambda_n\leq\Gamma_n\leq1\),
\[
\begin{aligned}
    \E[
        \widehat X_n^2
        \mid
        \mathscr F_{n-1}
    ]
    &\leq
    3
    \left[
        \|\nabla g(\theta_n)\|^2
        +a_\Gamma^2
        +\frac{\beta^2}{4}\zeta_n^2
    \right]
    \E[
        \Gamma_n
        \mid
        \mathscr F_{n-1}
    ].
\end{aligned}
\]
Indeed, the three terms are controlled respectively by
\(\|\nabla g(\theta_n)\|^2\E[\Gamma_n\mid\mathscr F_{n-1}]\),
\(a_\Gamma^2\E[\Gamma_n^4\mid\mathscr F_{n-1}]\), and
\((\beta^2/4)\zeta_n^2
\,\E[\Lambda_n^4\mid\mathscr F_{n-1}]\).
Now use
\[
    \|\nabla g(\theta_n)\|^2\leq2LV_n,
    \qquad
    \frac{\alpha_0\beta}{2}\zeta_n
    =V_n-g(\theta_n)\leq V_n,
\]
and multiply by \(\alpha_0^2\).  This proves
\cref{eq:xhat-conditional-variance}; the band-local bound follows
from \(V_n\leq2G\).
\end{namedproof}

\begin{namedproof}{Proof of \cref{lem:adj:ghat}.}
The AdaGrad-Norm update satisfies
\(\|\theta_{n+1}-\theta_n\|
=\alpha_0\|G_n\|/\sqrt{S_n}\leq\alpha_0\).
By smoothness and
\(\|\nabla g(\theta_n)\|^2\leq2LV_n\),
\[
    g(\theta_{n+1})-g(\theta_n)
    \leq
    \alpha_0\sqrt{2LV_n}
    +
    \frac{L\alpha_0^2}{2}.
\]
Moreover,
\[
\begin{aligned}
    &\|\nabla g(\theta_{n+1})\|^2
    -
    \|\nabla g(\theta_n)\|^2
     \leq
    2L\alpha_0\sqrt{2LV_n}
    +
    L^2\alpha_0^2.
\end{aligned}
\]
Because \(S_n\geq S_{n-1}\),
\[
\begin{aligned}
    V_{n+1}-V_n
    \leq{}&
    \alpha_0\sqrt{2LV_n}
    +
    \frac{L\alpha_0^2}{2}
    +
    \frac{\beta\alpha_0}{2\sqrt{S_0}}
    \left(
        2L\alpha_0\sqrt{2LV_n}
        +
        L^2\alpha_0^2
    \right).
\end{aligned}
\]
Hence \cref{eq:lyapunov-one-step-growth} holds with
\[
    a_h
    :=
    \sqrt{2L}
    \left(
        1+\frac{\beta\alpha_0L}{\sqrt{S_0}}
    \right)\alpha_0, \ \
    b_h
    :=
    \left(
        1+\frac{\beta\alpha_0L}{\sqrt{S_0}}
    \right)
    \frac{L\alpha_0^2}{2}.
\]
Since \(h(x)=O(\sqrt{x})\), there exists \(C_0<\infty\) such that
\(h(x)<x/2\) whenever \(x\geq C_0\).
\end{namedproof}

\begin{namedproof}{Proof of \cref{pro_0}.}
Choose \(0<\delta<
    \liminf_{\|\theta\|\to\infty}
    \|\nabla g(\theta)\|\).
Then $ K_\delta
    :=
    \{\theta:\|\nabla g(\theta)\|\leq\delta\}$
is compact. Consequently,
\[
    C_\delta
    :=
    \sup_{\theta\in K_\delta}g(\theta)
    +
    \frac{\beta\alpha_0\delta^2}{2\sqrt{S_0}}
    <
    \infty.
\]
Choose \(G_0>\max\{C_\delta,C_0\}\),
where \(C_0\) is given by
Lemma~\ref{lem:adj:ghat}. If \(V_n>G\geq G_0\), then
\(\theta_n\notin K_\delta\), proving part~(i).

For part~(ii), if \(V_{n-1}\leq G\), then the monotonicity of
\(h\) and \cref{eq:lyapunov-one-step-growth} give
\[
    V_n
    \leq
    V_{n-1}+h(V_{n-1})
    \leq
    G+h(G)
    <
    \frac32G.
\]
The same argument with \(2G\) in place of \(G\) proves part~(iii).
\end{namedproof}

\begin{namedproof}{Proof of \cref{lem:high-level-objective-variation}.}
The two-sided smoothness estimate and the update give
\[
    |g(\theta_{n+1})-g(\theta_n)|
    \leq
    \frac{\alpha_0\|\nabla g(\theta_n)\|\|G_n\|}
         {\sqrt{S_{n-1}}}
    +
    \frac{L\alpha_0^2}{2}
    \frac{\|G_n\|^2}{S_n}.
\]
On \(\{\|\nabla g(\theta_n)\|>\delta\}\), conditional
Cauchy--Schwarz and the affine second-moment bound imply $    \E[\|G_n\|\mid\mathscr F_{n-1}]
    \leq
    \sqrt{\beta+\frac{\sigma^2}{\delta^2}}\,
    \|\nabla g(\theta_n)\|.$
Moreover, $ 1/{S_n}
    \leq
    1/(\sqrt{S_0}\sqrt{S_{n-1}}),$
so the conditional expectation of the quadratic term is bounded by
\[
    \frac{L\alpha_0^2}{2\sqrt{S_0}}
    \left(
      \beta+\frac{\sigma^2}{\delta^2}
    \right)\zeta_n.
\]
Combining the two estimates proves
\cref{eq:adagrad-high-level-variation}.
\end{namedproof}

\section{Proofs of Preliminary Results}
\label{app:deferred-setup-proofs}

\subsection{Proof of Linear Growth from Non-Flatness}
\begin{namedproof}{Proof of \cref{lem:nonflat-linear-growth}.}
Choose \(R_c\) so that \(\|\nabla g(\theta)\|\geq c\) whenever
\(\|\theta\|\geq R_c\).  For \(\theta\notin B(0,R_c)\), let \(z(0)=\theta\)
and consider the globally Lipschitz vector field
\[
    \dot z(t)
    =
    -\frac1c
    \operatorname{proj}_{\overline B(0,c)}\!\bigl(\nabla g(z(t))\bigr).
\]
The solution exists for all \(t\geq0\) and satisfies
\(\|\dot z(t)\|\leq1\).  Until its first entrance time \(\tau\) into
\(\overline B(0,R_c)\), the projected vector has norm \(c\), and hence
\[
    \frac{\mathrm d}{\mathrm dt}g(z(t))
    =-\|\nabla g(z(t))\|
    \leq-c.
\]
The entrance time must be finite, since otherwise the preceding display
would contradict \(g\geq0\).  Moreover,
\(\tau\geq\operatorname{dist}(\theta,\overline B(0,R_c))
=\|\theta\|-R_c\).  Therefore
\[
    g(\theta)
    \geq g(z(\tau))+c\tau
    \geq c(\|\theta\|-R_c).
\]
For \(\theta\in\overline B(0,R_c)\), \cref{eq:nonflat-linear-growth}
follows from \(g\geq0\).
\end{namedproof}

\subsection{Proofs of the Local Relative-Moment Examples}
\begin{namedproof}{Proof of \cref{prop:ass:sto:bound}.}
Choose
\(0<\delta_0<
\liminf_{\|\theta\|\to\infty}\|\nabla g(\theta)\|\).
Then
\(K_{\delta_0}
=\{\theta\in\R^d:\|\nabla g(\theta)\|\leq\delta_0\}\)
is compact by
\cref{ass_g_poi}~\ref{ass_g_poi:i3}. For an admissible mini-batch \(\mathcal B\), define
\[
    G_{\mathcal B}(\theta)
    :=
    \frac1{|\mathcal B|}
    \sum_{i\in\mathcal B}\nabla g_i(\theta).
\]
Since every \(g_i\) is continuously differentiable,
\(G_{\mathcal B}\) is continuous and therefore bounded on
\(K_{\delta_0}\). Because the data set and the batch size are fixed, there are only
finitely many possible mini-batches, counting repeated indices when
sampling with replacement. Hence
\[
    \delta_1
    :=
    \max_{\mathcal B}
    \sup_{\theta\in K_{\delta_0}}
    \|G_{\mathcal B}(\theta)\|
    <
    \infty.
\]
It follows that
\(\|\nabla g(\theta)\|\leq\delta_0\) implies
\(\|G_{\mathcal B}(\theta)\|\leq\delta_1\) for every admissible
mini-batch \(\mathcal B\). Therefore
\(\|G_{\mathcal B}(\theta)\|^{2+\rho}
\leq\delta_1^\rho\|G_{\mathcal B}(\theta)\|^2\) on the local
region, proving \cref{ass:local-relative-moment} with
\(\kappa_\rho=\delta_1^\rho\).
\end{namedproof}
\begin{namedproof}{Proof of \cref{prop:additive-noise-lrm}.}
On \(\{\|\nabla g(\theta_n)\|\leq\delta_0\}\),
\[
 \E[\|G_n\|^{2+\rho}\mid\mathscr F_{n-1}]
 \leq 2^{1+\rho}(\delta_0^{2+\rho}+M_{2+\rho}),
 \qquad
 \E[\|G_n\|^2\mid\mathscr F_{n-1}]\geq v_0.
\]
The result follows with
\(\kappa_\rho=2^{1+\rho}
(\delta_0^{2+\rho}+M_{2+\rho})/v_0\).
\end{namedproof}

\section{Auxiliary Probabilistic Lemmas}
\label{sec:appendix-auxiliary}

We collect the auxiliary probabilistic tools used below.

\begin{lem}[Conditional summability]
\label{lem:conditional-summability}
Let \(\{Z_n\}_{n\geq1}\) be nonnegative and adapted, and set
\(q_n:=\E[Z_n\mid\mathscr F_{n-1}]\). Then, up to a null event,
\[
 \Bigl\{\sum_n q_n<\infty\Bigr\}
 \subseteq
 \Bigl\{\sum_n Z_n<\infty\Bigr\}.
\]
\end{lem}
\begin{proof}
For \(\ell\geq1\), let
\(\tau_\ell:=\inf\{N:\sum_{n=1}^Nq_n>\ell\}\). Predictability and Tonelli give
\(\E[\sum_{n<\tau_\ell}Z_n]=\E[\sum_{n<\tau_\ell}q_n]\leq\ell\), so the
stopped series is finite almost surely. On
\(\{\sum_nq_n\leq\ell\}\), \(\tau_\ell=\infty\); taking the union over
\(\ell\in\mathbb N\) proves the claim.
\end{proof}

The next lemma is a finite-dimensional martingale Rosenthal--Burkholder
inequality; see, for example, \citet{pinelis1994optimum}.

\begin{lem}[Martingale Rosenthal--Burkholder inequality]
\label{vital2}
Let \(\{Y_n\}_{n\geq1}\) be an \(\R^d\)-valued
martingale-difference sequence. For every \(p\geq2\), there is a constant
\(C_p<\infty\), depending only on \(p\), such that for every deterministic
\(N\),
\begin{align*}
\E\!\left[\max_{1\leq k\leq N}
    \left\|\sum_{n=1}^{k}Y_n\right\|^p\right]
&\leq C_p\E\!\left[
   \left(\sum_{n=1}^{N}
      \E[\|Y_n\|^2\mid\mathscr F_{n-1}]
   \right)^{p/2}\right]+C_p\E\!\left[\sum_{n=1}^{N}
      \E[\|Y_n\|^p\mid\mathscr F_{n-1}]\right].
\end{align*}
The same estimate applies to a finite predictable truncation of a
martingale-difference sequence.
\end{lem}

\section{Lyapunov and Self-Normalization Estimates}
Throughout this section we retain the notation
\(G_n,\zeta_n,\Gamma_n,\widetilde\Gamma_n,V_n\) from
\cref{eq:adagrad-normalized-quantities,eq:adagrad-lyapunov}.
\subsection{Proof of the Sufficient Descent Lemma}
\label{proof:sufficient:lem}

\begin{namedproof}{Proof of \cref{sufficient:lem}.}
By \(L\)-smoothness,
\begin{equation}
\label{eq:smoothness-descent-proof}
    g(\theta_{n+1})-g(\theta_n)
    \leq
    -\alpha_0
    \frac{\nabla g(\theta_n)^\top G_n}{\sqrt{S_n}}
    +
    \frac{L\alpha_0^2}{2}\Gamma_n.
\end{equation}
Since \(S_n\) depends on the current oracle output, introduce
\[
    \Lambda_n
    :=
    \frac{\|G_n\|^2}
    {\sqrt{S_n}
     \bigl(\sqrt{S_{n-1}}+\sqrt{S_n}\bigr)}.
\]
Then
\[
    \frac{1}{\sqrt{S_{n-1}}}
    -
    \frac{1}{\sqrt{S_n}}
    =
    \frac{\Lambda_n}{\sqrt{S_{n-1}}},
    \qquad
    0\leq\Lambda_n\leq\Gamma_n\leq1.
\]

Adding and subtracting the predictable-denominator term and using
conditional unbiasedness gives
\begin{align}
g(\theta_{n+1})-g(\theta_n)
\leq{}&
-\alpha_0\zeta_n
+
\alpha_0
\E\left[
    \frac{\|\nabla g(\theta_n)\|\|G_n\|}
         {\sqrt{S_{n-1}}}
    \Lambda_n
    \,\middle|\,
    \mathscr F_{n-1}
\right]
+
\frac{L\alpha_0^2}{2}\Gamma_n
+
\alpha_0X_n,
\label{eq:descent-predictable-decomposition}
\end{align}
where
\[
    X_n
    :=
    \E\left[
        \frac{\nabla g(\theta_n)^\top G_n}{\sqrt{S_n}}
        \,\middle|\,
        \mathscr F_{n-1}
    \right]
    -
    \frac{\nabla g(\theta_n)^\top G_n}{\sqrt{S_n}}
\]
is a martingale difference.

Conditional Cauchy--Schwarz, Young's inequality, and the affine
second-moment assumption yield
\begin{align}
&\E\left[
    \frac{\|\nabla g(\theta_n)\|\|G_n\|}
         {\sqrt{S_{n-1}}}
    \Lambda_n
    \,\middle|\,
    \mathscr F_{n-1}
\right] \leq
\frac12\zeta_n
+
\frac{\sigma^2}{2\sqrt{S_0}}
\E[\Gamma_n^2\mid\mathscr F_{n-1}]
+
\frac{\beta}{2}\zeta_n
\E[\Lambda_n^2\mid\mathscr F_{n-1}].
\label{eq:adaptive-correction-bound}
\end{align}
Let \(a_\Gamma=\sigma^2/(2\sqrt{S_0})\), and define explicitly
\[
\begin{aligned}
    \widehat X_n
    :={}&
    X_n
    +
    a_\Gamma
    \left(
      \E[\Gamma_n^2\mid\mathscr F_{n-1}]-\Gamma_n^2
    \right)
    +
    \frac{\beta}{2}\zeta_n
    \left(
      \E[\Lambda_n^2\mid\mathscr F_{n-1}]-\Lambda_n^2
    \right).
\end{aligned}
\]
Every added term is conditionally centered.  Substituting this
identity into \cref{eq:adaptive-correction-bound} gives
\begin{align}
g(\theta_{n+1})-g(\theta_n)
\leq{}&
-\frac{\alpha_0}{2}\zeta_n
+
\left(
    \frac{\alpha_0\sigma^2}{2\sqrt{S_0}}
    +
    \frac{L\alpha_0^2}{2}
\right)\Gamma_n
+
\frac{\beta\alpha_0}{2}
\zeta_n\Lambda_n^2
+
\alpha_0\widehat X_n,
\label{eq:pre-lyapunov-descent}
\end{align}
where \(\Gamma_n^2\leq\Gamma_n\) was used.

Next,
\begin{align}
    \zeta_n\Lambda_n^2
    &\leq
    \|\nabla g(\theta_n)\|^2
    \left(
        \frac{1}{\sqrt{S_{n-1}}}
        -
        \frac{1}{\sqrt{S_n}}
    \right)
    =
    \zeta_n-\zeta_{n+1}
    +
    \frac{
        \|\nabla g(\theta_{n+1})\|^2
        -
        \|\nabla g(\theta_n)\|^2
    }{\sqrt{S_n}}.
\label{eq:zeta-lambda-telescope}
\end{align}
For every \(\varepsilon>0\), smoothness and Young's inequality imply
\begin{equation}
\label{eq:gradient-difference-general}
    \frac{
        \|\nabla g(\theta_{n+1})\|^2
        -
        \|\nabla g(\theta_n)\|^2
    }{\sqrt{S_n}}
    \leq
    \varepsilon\zeta_n
    +
    L^2\alpha_0^2
    \left(1+\frac1\varepsilon\right)
    \widetilde\Gamma_n.
\end{equation}
Choose $ \varepsilon
    :=
    \frac{1}{2(1+\beta)}.$
Then $ \frac{\beta\varepsilon}{2}\leq\frac14.$
Define $ M_n:=\alpha_0\widehat X_n.$
Combining the preceding estimates gives
\begin{equation}
\label{eq:sufficient-descent-final}
    V_{n+1}-V_n
    \leq
    -\frac{\alpha_0}{4}\zeta_n
    +
    C_{\Gamma,1}\Gamma_n
    +
    C_{\Gamma,2}\widetilde\Gamma_n
    +
    M_n,
\end{equation}
where $ C_{\Gamma,1}
    :=
    \frac{\alpha_0\sigma^2}{2\sqrt{S_0}}
    +
    \frac{L\alpha_0^2}{2},$
and \(C_{\Gamma,2}<\infty\) depends only on
\(\alpha_0,L,\beta\). Moreover, $ \E[M_n\mid\mathscr F_{n-1}]=0.$
\end{namedproof}

\subsection{Proof of the Normalized Gradient-Budget Estimate}
\label{sec:proof:lemsu}

\begin{namedproof}{Proof of \cref{lem_su}.}
Let $   y_n:=\frac{1}{\sqrt{S_{n-1}}}.$
Multiplying \cref{eq:sufficient-descent-final} by \(y_n\), summing, and
taking expectations gives
\begin{align}
\frac{\alpha_0}{4}
\E\left[
    \sum_{n=1}^{T}y_n\zeta_n
\right]
\leq{}&
y_1V_1
+
C_{\Gamma,1}
\E\left[
    \sum_{n=1}^{T}y_n\Gamma_n
\right]
+
C_{\Gamma,2}
\E\left[
    \sum_{n=1}^{T}y_n\widetilde\Gamma_n
\right].
\label{eq:standard-zeta-summability}
\end{align}
The martingale term vanishes because \(y_n\) is predictable.  The
summation-by-parts identity
\[
    \sum_{n=1}^{T}y_n(V_n-V_{n+1})
    =y_1V_1-y_TV_{T+1}
     +\sum_{n=2}^{T}(y_n-y_{n-1})V_n
    \leq y_1V_1
\]
uses \(V_n\geq0\) and the monotonicity \(y_n\leq y_{n-1}\), and
justifies the first term on the right-hand side of
\cref{eq:standard-zeta-summability}.

The remaining self-normalized sums are deterministic telescoping budgets.
Indeed, with \(y_{n+1}=S_n^{-1/2}\),
\[
    y_n\Gamma_n
    =\frac{S_n-S_{n-1}}{\sqrt{S_{n-1}}S_n}
    \leq2(y_n-y_{n+1}),
\]
and
\[
    y_n\widetilde\Gamma_n
    =\frac{S_n-S_{n-1}}{\sqrt{S_{n-1}}S_n^{3/2}}
    \leq
    \frac1{S_{n-1}}-\frac1{S_n}.
\]
Consequently,
\[
    \sum_{n=1}^{\infty}y_n\Gamma_n
    \leq
    \frac{2}{\sqrt{S_0}},
    \qquad
    \sum_{n=1}^{\infty}y_n\widetilde\Gamma_n
    \leq
    \frac{1}{S_0}.
\]
Therefore, for some \(C_\zeta<\infty\),
\begin{equation}
\label{eq:gradient-energy-standard}
    \E\left[
        \sum_{n=1}^{\infty}
        \frac{\|\nabla g(\theta_n)\|^2}{S_{n-1}}
    \right]
    \leq
    C_\zeta.
\end{equation}
Since the series is nonnegative, it is also finite almost surely.

For every \(\delta>0\), conditional unbiasedness and the affine
second-moment bound imply
\begin{align}
&\E\left[
    \sum_{n=1}^{\infty}
    \ind_{\{\|\nabla g(\theta_n)\|>\delta\}}
    \frac{\|G_n\|^2}{S_{n-1}}
\right]\leq
\left(
    \beta+\frac{\sigma^2}{\delta^2}
\right)
\E\left[
    \sum_{n=1}^{\infty}
    \frac{\|\nabla g(\theta_n)\|^2}{S_{n-1}}
\right]
<\infty.
\label{eq:high-gradient-oracle-energy}
\end{align}

Since \(S_n\geq S_{n-1}\),
\[
    \Gamma_n
    \leq
    \frac{\|G_n\|^2}{S_{n-1}},
\]
and hence
\begin{equation}
\label{eq:high-gradient-gamma-energy}
    \sum_{n=1}^{\infty}
    \E\left[
        \ind_{\{\|\nabla g(\theta_n)\|>\delta\}}
        \Gamma_n
    \right]
    \leq
    \left(\beta+\frac{\sigma^2}{\delta^2}\right)C_\zeta
    <\infty.
\end{equation}
Thus \cref{eq:adagrad-high-gradient-energy} holds, for example, with
\(M_\delta=C_\zeta\).
\end{namedproof}

\section{ODE and Asymptotic-Pseudotrajectory Tools}
\label{app:ode-identification}

The following standard consequence of the asymptotic-pseudotrajectory framework
\citep{benaim2006dynamics} is the only ODE input used in the main text.

\begin{pros}[ODE identification principle]
\label{SA_p}
Consider the recursion \cref{SA}. Suppose that the following conditions hold
almost surely.

\begin{enumerate}[label=\textnormal{(A.\arabic*)},leftmargin=*]

    \item
    \label{pros:a1}
    The sequence \(\{z_n\}_{n\geq1}\) is bounded.

    \item
    \label{pros:a2}
    The stepsizes satisfy
    \[
        \gamma_n\longrightarrow0,
        \qquad
        \sum_{n=1}^{\infty}\gamma_n=\infty,
    \]
    and, for every \(T>0\),
    \begin{equation}
    \label{eq:finite-window-perturbation}
        \lim_{n\to\infty}
        \sup_{n\leq k\leq m(t_n+T)}
        \left\|
            \sum_{i=n}^{k}\gamma_iU_i
        \right\|
        =0.
    \end{equation}

    \item
    \label{pros:a3}
    The limiting ODE admits a strict Lyapunov function whose values on the
    equilibrium set
    \[
        \mathcal E
        :=
        \{z\in\R^d:F(z)=0\}
    \]
    have empty interior.

\end{enumerate}

Then the limit set of \(\{z_n\}\) is contained in \(\mathcal E\).  In
particular, every accumulation point of \(\{z_n\}\) is an equilibrium of the
limiting ODE.
\end{pros}

The proposition identifies accumulation points rather than a single limit;
in the gradient-flow application, boundedness and continuity of \(\nabla g\)
are enough to deduce convergence of the full gradient-norm sequence.

\section{Auxiliary Results for Asymptotic Stationarity}
\label{sec:conv1:aux}
\begin{lem}\label{lem:effective-time-divergence}
If \cref{ass_g_poi,ass_noise} hold and \(\gamma_n:=\alpha_0/\sqrt{S_n}\), then $\sum_{n=1}^{+\infty}\gamma_n=+\infty \quad \text{a.s.}$
\end{lem}

\subsection{Divergence of the Effective Time}
\label{proof:step-size}

\begin{namedproof}{Proof of \cref{lem:effective-time-divergence}.}
Recall that $  \gamma_n=\frac{\alpha_0}{\sqrt{S_n}}.$
If \(S_n\) remains bounded, then $\sum_{n=1}^{\infty}\gamma_n=\infty$
immediately.  It therefore remains to rule out finite effective time on the
event \(\{S_n\to\infty\}\).

Define the bad event
\[
    \mathcal B_{\mathrm{eff}}
    :=
    \left\{
      S_n\to\infty,
      \ \sum_{n=1}^{\infty}\frac1{\sqrt{S_n}}<\infty
    \right\}.
\]
On \(\mathcal B_{\mathrm{eff}}\), \cref{eq:smooth-gradient-objective-bound} together with
\cref{thm:main}~\ref{thm:main:stability} give
\[
\begin{aligned}
    \sum_{n=1}^{\infty}
    \frac{\|\nabla g(\theta_{n+1})\|^2}{\sqrt{S_n}}
    &\leq
    2L\sup_{j\geq1}g(\theta_j)
    \sum_{n=1}^{\infty}\frac1{\sqrt{S_n}}
    <\infty
    \qquad\text{almost surely on }\mathcal B_{\mathrm{eff}}.
\end{aligned}
\]
The affine conditional second-moment assumption therefore gives
\[
    \sum_{n=1}^{\infty}
    \E\left[
        \frac{\|G_{n+1}\|^2}{\sqrt{S_n}}
        \,\middle|\,
        \mathscr F_n
    \right]
    <
    \infty
    \qquad\text{almost surely on }\mathcal B_{\mathrm{eff}}.
\]
Applying \cref{lem:conditional-summability} after the one-step
reindexing
\(\widehat Z_1:=0\) and
\(\widehat Z_{n+1}:=\|G_{n+1}\|^2/\sqrt{S_n}\), \(n\geq1\), gives
\[
    \sum_{n=1}^{\infty}
    \frac{\|G_{n+1}\|^2}{\sqrt{S_n}}
    <
    \infty
    \qquad\text{almost surely on }\mathcal B_{\mathrm{eff}}.
\]

On the other hand,
\[
\begin{aligned}
    \frac{\|G_{n+1}\|^2}{\sqrt{S_n}}
    &=
    \frac{S_{n+1}-S_n}{\sqrt{S_n}}
    \geq
    \int_{S_n}^{S_{n+1}}x^{-1/2}\,\mathrm dx
    =
    2\bigl(\sqrt{S_{n+1}}-\sqrt{S_n}\bigr).
\end{aligned}
\]
Summing over \(n\) shows that the left-hand series diverges whenever
\(S_n\to\infty\).  Hence \(\mathbb P(\mathcal B_{\mathrm{eff}})=0\),
and therefore $\sum_{n=1}^{\infty}\gamma_n=\infty
    \qquad\text{almost surely}.$
\end{namedproof}

\section{Relative-Moment Estimates and Finite-Window Noise Control}
\label{app:local-relative-moment-proofs}

\subsection{Why an Ordinary Local Higher-Moment Bound Is Insufficient}
\label{app:ordinary-local-moment-counterexample}

\begin{namedproof}{Proof of \cref{prop:ordinary-local-moment-insufficient}.}
Fix \(p=2+\rho>2\) and let \(g(x)=x^2/2\), \(\theta_1=0\). The objective is
\(1\)-smooth, nonnegative, coercive, non-flat at infinity, and has critical
value set \(\{0\}\). Fix \(\alpha_0>0\), choose \(\eta\in(0,1)\), set
\(c_\eta:=(1-\eta)/\sqrt{1+(1+\eta)^2}\), and choose
\begin{equation}
\label{eq:counterexample-parameters}
 0<\varepsilon<\frac{\alpha_0c_\eta}{2},
 \qquad
 S_0>\max\left\{1,\alpha_0^2,\frac{\varepsilon^2}{\eta^2}\right\}.
\end{equation}
Set \(\delta_0:=\varepsilon\). Given \(\mathscr F_{n-1}\), put
\(\nu_n=0\) when \(|\theta_n|>\varepsilon\); otherwise let
\begin{equation}
\label{eq:counterexample-oracle-law}
\nu_n=
\begin{cases}
+\sqrt{S_{n-1}},&\text{with probability }(2S_{n-1}^{p/2})^{-1},\\
-\sqrt{S_{n-1}},&\text{with probability }(2S_{n-1}^{p/2})^{-1},\\
0,&\text{with probability }1-S_{n-1}^{-p/2},
\end{cases}
\end{equation}
Let
\(G_n=\theta_n+\nu_n\), \(S_n=S_{n-1}+G_n^2\), and
\(\theta_{n+1}=\theta_n-\alpha_0G_n/\sqrt{S_n}\).

\paragraph{Oracle moments.}
Symmetry gives \(\E[G_n\mid\mathscr F_{n-1}]=\theta_n\), while
\[
 \E[G_n^2\mid\mathscr F_{n-1}]
 =\theta_n^2+
 \ind_{\{|\theta_n|\leq\varepsilon\}}S_{n-1}^{1-p/2}
 \leq\theta_n^2+S_0^{-\rho/2}.
\]
Thus the affine second-moment condition holds. On the local region,
\(\E[|\nu_n|^p\mid\mathscr F_{n-1}]=1\), hence
\(
 \E[|G_n|^p\mid\mathscr F_{n-1}]
 \leq2^{p-1}(\varepsilon^p+1)
\); the ordinary local \(p\)-moment is uniformly bounded.

\paragraph{Spike-free contraction.}
Without a spike,
\(
S_n=S_{n-1}+\theta_n^2
\)
and
\(
\theta_{n+1}=\theta_n(1-\alpha_0/\sqrt{S_n})
\).
Since \(S_0>\alpha_0^2\), the factor lies in \((0,1)\). Moreover,
\[
 \theta_n^2-\theta_{n+1}^2
 \geq\frac{\alpha_0\theta_n^2}{\sqrt{S_n}},
 \qquad
 \sqrt{S_n}-\sqrt{S_{n-1}}
 \leq\frac{\theta_n^2}{\sqrt{S_n}}.
\]
Hence on any spike-free interval beginning at \(m\),
\begin{equation}
\label{eq:counterexample-spike-free-budget}
 \sqrt{S_k}-\sqrt{S_{m-1}}\leq\frac{\theta_m^2}{\alpha_0},
 \qquad k\geq m.
\end{equation}
The accumulator therefore stays bounded on such a continuation, so
\(|\theta_k|\) contracts geometrically, enters
\([-\varepsilon,\varepsilon]\) in finite time, and stays there until the next
spike.

\paragraph{Infinitely many spikes.}
Let \(T_0=0\) and
\(T_{j+1}:=\inf\{n>T_j:\nu_n\neq0\}\). On \(\{T_j<\infty\}\), condition on
\(\mathscr F_{T_j}\) and consider a continuation with no further spike.
By \cref{eq:counterexample-spike-free-budget}, its accumulator is bounded by
some finite random \(\overline S_j\); after re-entering the local region, the
conditional spike probability is therefore at least
\(c_j:=\overline S_j^{-p/2}>0\). Consequently,
\[
 \Pro(T_{j+1}=\infty\mid\mathscr F_{T_j})
 \leq\lim_{N\to\infty}(1-c_j)^N=0,
\]
and induction shows that spikes occur infinitely often almost surely.

\paragraph{Failure of stationarity.}
At a spike write \(y_n:=\theta_n/\sqrt{S_{n-1}}\) and
\(s_n\in\{-1,+1\}\) for its sign. By
\cref{eq:counterexample-parameters}, \(|y_n|\leq\eta\), and
\[
 \theta_{n+1}
 =\theta_n-\alpha_0\frac{y_n+s_n}{\sqrt{1+(y_n+s_n)^2}},
 \qquad
 |\theta_{n+1}|\geq\alpha_0c_\eta-\varepsilon
 \geq\frac{\alpha_0c_\eta}{2}>0.
\]
Thus \(|\nabla g(\theta_n)|\not\to0\). Every normalized update has magnitude
at most \(\alpha_0\), and between spikes \(|\theta_n|\) decreases, so the
iterates remain bounded. Finally, each spike satisfies
\(
 S_n\geq[1+(1-\eta)^2]S_{n-1}
\), hence \(S_n\to\infty\). This proves the proposition.
\end{namedproof}

\begin{rem}[Scope of the counterexample]
The law in \cref{eq:counterexample-oracle-law} may depend on \(S_{n-1}\), as
allowed by \cref{ass_noise}. The example therefore concerns the general
adapted-oracle model and makes no claim about fresh i.i.d. state-only oracles.
\end{rem}

\subsection{Weighted Local Second-Moment Control}
\label{app:weighted-local-second-moment}

For the remainder of this section, set
\begin{equation}
\label{eq:lrm-notation}
    h_n:=\nabla g(\theta_n),
    \qquad
    \varepsilon_n:=G_n-h_n,
    \qquad
    L_n:=\ind_{\{\|h_n\|\leq\delta_0\}},
    \qquad
    \bar L_n:=1-L_n,
\end{equation}
and define the predictable local conditional second moment
\begin{equation}
\label{eq:lrm-local-energy}
    v_n:=L_n\E[\|G_n\|^2\mid\mathscr F_{n-1}].
\end{equation}

\begin{lem}[Weighted local second-moment control]
\label{lem:weighted-local-second-moment}
Suppose that conditional unbiasedness, the affine conditional second-moment
bound, and \cref{ass:local-relative-moment} hold. Put
$\bar\rho:=\min\{\rho,1\}$ and $p:=2+\bar\rho$. There exists a finite constant
$K_p$ such that, almost surely for every $n$,
\begin{equation}
\label{eq:lrm-relative-p-moment}
    L_n\E[\|G_n\|^p\mid\mathscr F_{n-1}]
    +L_n\E[\|\varepsilon_n\|^p\mid\mathscr F_{n-1}]
    \leq K_pv_n.
\end{equation}
Moreover, for every $a>0$,
\begin{equation}
\label{eq:lrm-weighted-energy}
    \E\left[
      \sum_{n=1}^{\infty}\frac{v_n}{S_{n-1}^{1+a}}
    \right]<\infty,
\end{equation}
and consequently
\begin{equation}
\label{eq:lrm-normalized-spikes}
    \sum_{n=1}^{\infty}
    L_n\left(\frac{\|G_n\|}{\sqrt{S_{n-1}}}\right)^p
    <\infty
    \qquad\text{almost surely}.
\end{equation}
\end{lem}

\begin{proof}
If $\rho\leq1$, the first $p$-moment term in
\cref{eq:lrm-relative-p-moment} follows directly from
\cref{ass:local-relative-moment}. If $\rho>1$, use
$x^3\leq x^2+x^{2+\rho}$ for $x\geq0$. Conditional unbiasedness and
Jensen's inequality give
$\|h_n\|^2\leq\E[\|G_n\|^2\mid\mathscr F_{n-1}]$. On $\{L_n=1\}$,
$\|h_n\|\leq\delta_0$, so
\[
\begin{aligned}
L_n\E[\|\varepsilon_n\|^p\mid\mathscr F_{n-1}]
&\leq
2^{p-1}L_n\left(
  \E[\|G_n\|^p\mid\mathscr F_{n-1}]+\|h_n\|^p
\right)\\
&\leq K_pv_n
\end{aligned}
\]
for a finite $K_p$, after enlarging the constant if necessary.  This proves
\cref{eq:lrm-relative-p-moment}.

Fix $a>0$, $N\geq1$, and $R>1$, and set
\[
    A_{a,N}
    :=
    \E\left[
      \sum_{n=1}^{N}\frac{L_n\|G_n\|^2}{S_{n-1}^{1+a}}
    \right]
    =
    \E\left[
      \sum_{n=1}^{N}\frac{v_n}{S_{n-1}^{1+a}}
    \right].
\]
On $\{\|G_n\|\leq R\sqrt{S_{n-1}}\}$, one has
$S_n\leq(1+R^2)S_{n-1}$, and hence
\[
    L_n\ind_{\{\|G_n\|\leq R\sqrt{S_{n-1}}\}}
    \frac{\|G_n\|^2}{S_{n-1}^{1+a}}
    \leq
    (1+R^2)^{1+a}\frac{S_n-S_{n-1}}{S_n^{1+a}}.
\]
Summing the right-hand side and comparing with an integral gives
\[
    \sum_n
    (1+R^2)^{1+a}\frac{S_n-S_{n-1}}{S_n^{1+a}}
    \leq
    \frac{(1+R^2)^{1+a}}{aS_0^a}.
\]
On the complementary event,
\[
    L_n\ind_{\{\|G_n\|>R\sqrt{S_{n-1}}\}}
    \frac{\|G_n\|^2}{S_{n-1}^{1+a}}
    \leq
    R^{-\rho}
    L_n\frac{\|G_n\|^{2+\rho}}
         {S_{n-1}^{1+a+\rho/2}}.
\]
Taking conditional expectations and using
\cref{ass:local-relative-moment} yields
\[
    A_{a,N}
    \leq
    \frac{(1+R^2)^{1+a}}{aS_0^a}
    +\kappa_\rho R^{-\rho}S_0^{-\rho/2}A_{a,N}.
\]
Choose $R$ sufficiently large that the final coefficient is at most $1/2$.
The resulting estimate is uniform in $N$, and monotone convergence gives
\cref{eq:lrm-weighted-energy}.

Finally, \cref{eq:lrm-relative-p-moment,eq:lrm-weighted-energy}, with
$a=\bar\rho/2$, imply $   \E\left[
      \sum_{n=1}^{\infty}
      L_n\frac{\|G_n\|^p}{S_{n-1}^{p/2}}
    \right]<\infty.$
The nonnegative series is therefore finite almost surely, proving
\cref{eq:lrm-normalized-spikes}.
\end{proof}

\subsection{Finite-Window Perturbations Without Bounded Oracle Values}
\label{app:lrm-finite-window}

\begin{lem}[Finite-window perturbations under the local relative-moment condition]
\label{inequ:gammaU}
Assume \cref{ass_g_poi,ass_noise,ass:local-relative-moment}. On
$\mathcal A^c:=\{S_n\to\infty\}$, set
\(\gamma_n:=\alpha_0/\sqrt{S_n}\),
\(t_0:=0\),
\(t_n:=\sum_{i=1}^n\gamma_i\), and
\(m(t):=\max\{j\ge0:t_j\le t\}\).
Then, for every $T>0$,
\begin{equation}
\label{eq:lrm-finite-window}
    \lim_{n\to\infty}
    \sup_{n\leq k\leq m(t_n+T)}
    \left\|
      \sum_{i=n}^{k}\gamma_i\varepsilon_i
    \right\|
    =0
    \qquad\text{almost surely on }\mathcal A^c,
\end{equation}
where $\varepsilon_i=G_i-\nabla g(\theta_i)$.
\end{lem}

\begin{proof}
The perturbation in \cref{eq:adagrad-direct-sa} is
$-\varepsilon_n$, and therefore its sign is immaterial for
\cref{eq:lrm-finite-window}.

\paragraph{Step 1: relative jumps and the denominator correction.}
Put
$r_n:=\|G_n\|/\sqrt{S_{n-1}}$. By
\cref{eq:lrm-normalized-spikes}, $L_nr_n\to0$ almost surely. Taking
$\delta=\delta_0$ in \cref{eq:high-gradient-oracle-energy} gives
\begin{equation}
\label{eq:lrm-high-gradient-energy}
    \E\left[
      \sum_{n=1}^{\infty}\bar L_n\frac{\|G_n\|^2}{S_{n-1}}
    \right]<\infty.
\end{equation}
Hence $\bar L_nr_n\to0$ almost surely, and therefore, on $\mathcal A^c$,
\begin{equation}
\label{eq:lrm-relative-jump}
    r_n\to0,
    \qquad
    \sqrt{\frac{S_{n-1}}{S_n}}
    =(1+r_n^2)^{-1/2}\to1.
\end{equation}
For every finite interval of indices,
\begin{equation}
\label{eq:lrm-denominator-decomposition}
    \sum_i\frac{\varepsilon_i}{\sqrt{S_i}}
    =
    \sum_i\frac{\varepsilon_i}{\sqrt{S_{i-1}}}
    -
    \sum_i\left(
       \frac1{\sqrt{S_{i-1}}}-\frac1{\sqrt{S_i}}
    \right)\varepsilon_i.
\end{equation}
The denominator-correction series is absolutely convergent. On the local region,
\[
\begin{aligned}
&L_n\left(
  \frac1{\sqrt{S_{n-1}}}-\frac1{\sqrt{S_n}}
\right)\|\varepsilon_n\|\leq
\delta_0\left(
  \frac1{\sqrt{S_{n-1}}}-\frac1{\sqrt{S_n}}
\right)
+L_n\phi(r_n),
\end{aligned}
\]
where $\phi(r):=r\bigl(1-(1+r^2)^{-1/2}\bigr)$. For
$p=2+\bar\rho\in(2,3]$, there is $C_p<\infty$ such that
$\phi(r)\leq C_pr^p$ for all $r\geq0$: use
$\phi(r)\leq r^3/2$ for $r\leq1$ and $\phi(r)\leq r$ for $r\geq1$.
The first series telescopes and the second is finite by
\cref{eq:lrm-normalized-spikes}.

On the high-gradient region,
\[
    \bar L_n\left(
      \frac1{\sqrt{S_{n-1}}}-\frac1{\sqrt{S_n}}
    \right)\|\varepsilon_n\|
    \leq
    \bar L_n\frac{\|G_n\|\,\|\varepsilon_n\|}{S_{n-1}}.
\]
Conditional unbiasedness gives $\E[\|\varepsilon_n\|^2\mid\mathscr F_{n-1}]
    =
    \E[\|G_n\|^2\mid\mathscr F_{n-1}]-\|h_n\|^2,$
and conditional Cauchy--Schwarz therefore yields $ \E[\|G_n\|\,\|\varepsilon_n\|
       \mid\mathscr F_{n-1}]
    \leq
    \E[\|G_n\|^2\mid\mathscr F_{n-1}].$
Together with \cref{eq:lrm-high-gradient-energy} and Tonelli's theorem,
this proves absolute summability of the high-gradient correction. Thus, with
\[
    C_n
    :=
    \left(
      \frac1{\sqrt{S_{n-1}}}-\frac1{\sqrt{S_n}}
    \right)\varepsilon_n,
\]
we have
\begin{equation}
\label{eq:lrm-summable-correction}
    \sum_{n=1}^{\infty}\|C_n\|<\infty
    \quad\text{a.s.},
    \qquad
    \sup_{k\geq n}\left\|\sum_{i=n}^{k}C_i\right\|\to0.
\end{equation}

\paragraph{Step 2: first-passage blocks and the high-gradient martingale.}
Fix \(T>0\) and define
\(
 \tau_\ell:=\inf\{n\geq0:t_n\geq\ell T\}
\), \(\ell\geq0\). Each \(\tau_\ell\) is a stopping time, and the block
selector \(\ind_{\{\tau_\ell<n\leq\tau_{\ell+1}\}}\) is predictable. By
\cref{lem:effective-time-divergence}, every level is eventually reached. For
each nonempty block,
\begin{equation}
\label{eq:lrm-block-time}
 \sum_{n=\tau_\ell+1}^{\tau_{\ell+1}}\frac1{\sqrt{S_{n-1}}}
 \leq C_T,
 \qquad C_T:=\frac{T}{\alpha_0}+\frac{2}{\sqrt{S_0}}.
\end{equation}
Indeed, the analogous sum with \(S_n\) is at most
\(T/\alpha_0+S_0^{-1/2}\), and the shift from \(S_n\) to \(S_{n-1}\)
telescopes.

Define the high-gradient block maximum
\[
 B_{\ell,H}:=
 \sup_{\tau_\ell<k\leq\tau_{\ell+1}}
 \left\|\sum_{n=\tau_\ell+1}^{k}
 \bar L_n\frac{\varepsilon_n}{\sqrt{S_{n-1}}}\right\|,
\]
with value zero for empty blocks. Applying Doob's \(L^2\) maximal inequality
to deterministic truncations and then letting the truncation tend to infinity
gives
\begin{equation}
\label{eq:lrm-high-block-maxima}
 \sum_{\ell=0}^{\infty}\E[B_{\ell,H}^2]
 \leq4\E\sum_{n=1}^{\infty}
 \bar L_n\frac{\E[\|\varepsilon_n\|^2\mid\mathscr F_{n-1}]}{S_{n-1}}
 \leq4\E\sum_{n=1}^{\infty}
 \bar L_n\frac{\|G_n\|^2}{S_{n-1}}<\infty.
\end{equation}
The last inequality uses conditional unbiasedness and
\cref{eq:lrm-high-gradient-energy}. Hence \(B_{\ell,H}\to0\) almost surely.

\paragraph{Step 3: dyadic tiles for the low-gradient martingale.}
Let \(b_j:=2^jS_0\) and
\(\mathcal B_j:=\{n:b_j\leq S_{n-1}<2b_j\}\). For a block--band tile define
\[
 X_{\ell,j}:=
 \sup_{\tau_\ell<k\leq\tau_{\ell+1}}
 \left\|\sum_{\substack{\tau_\ell<n\leq k\\ n\in\mathcal B_j}}
 L_n\frac{\varepsilon_n}{\sqrt{S_{n-1}}}\right\|,
 \qquad
 V_{\ell,j}:=
 \sum_{\substack{\tau_\ell<n\leq\tau_{\ell+1}\\ n\in\mathcal B_j}}
 L_n\frac{\E[\|\varepsilon_n\|^2\mid\mathscr F_{n-1}]}{S_{n-1}}.
\]
On the local region, \(v_n\leq M_0:=\beta\delta_0^2+\sigma^2\), so
\cref{eq:lrm-block-time} yields
\begin{equation}
\label{eq:lrm-tile-quadratic-variation}
 V_{\ell,j}\leq M_0C_Tb_j^{-1/2}.
\end{equation}
Applying the martingale Rosenthal--Burkholder inequality in
Lemma~\ref{vital2} first to deterministic truncations and then passing
to the limit, with \(p=2+\bar\rho\), gives 
\begin{align}
\E[X_{\ell,j}^p]
&\leq C_p\E\left[V_{\ell,j}^{p/2}
 +\sum_{\substack{\tau_\ell<n\leq\tau_{\ell+1}\\ n\in\mathcal B_j}}
 L_n\frac{\E[\|\varepsilon_n\|^p\mid\mathscr F_{n-1}]}{S_{n-1}^{p/2}}
 \right] \notag\\
&\leq C_{p,T}\E\sum_{\substack{\tau_\ell<n\leq\tau_{\ell+1}\\ n\in\mathcal B_j}}
 v_n\left(S_{n-1}^{-1-\bar\rho/4}
          +S_{n-1}^{-1-\bar\rho/2}\right).
\label{eq:lrm-tile-moment}
\end{align}
Here the first term uses
\(V_{\ell,j}^{p/2}\leq C_{p,T}b_j^{-\bar\rho/4}V_{\ell,j}\), and
\cref{eq:lrm-relative-p-moment} controls the second.

Because the nondecreasing accumulator visits dyadic bands in increasing
order, the low-gradient block maximum
\[
 B_{\ell,L}:=
 \sup_{\tau_\ell<k\leq\tau_{\ell+1}}
 \left\|\sum_{n=\tau_\ell+1}^{k}
 L_n\frac{\varepsilon_n}{\sqrt{S_{n-1}}}\right\|
\]
satisfies \(B_{\ell,L}\leq\sum_{j\geq j_0(\ell)}X_{\ell,j}\), where
\(j_0(\ell)\) is the first visited band. For any
\(0<\eta<\bar\rho/4\), weighted H\"older gives
\[
 B_{\ell,L}^p
 \leq C_{p,\eta}\sum_{j\geq j_0(\ell)}
 2^{\eta(j-j_0(\ell))}X_{\ell,j}^p.
\]
Using the disjoint block--band partition together with
\cref{eq:lrm-tile-moment,eq:lrm-weighted-energy} yields
\begin{equation}
\label{eq:lrm-low-block-maxima}
 \sum_{\ell=0}^{\infty}\E[B_{\ell,L}^p]
 \leq C\E\sum_{n=1}^{\infty}v_n\left(
 S_{n-1}^{-1-(\bar\rho/4-\eta)}
 +S_{n-1}^{-1-(\bar\rho/2-\eta)}\right)<\infty.
\end{equation}
Thus \(B_{\ell,L}\to0\) almost surely. Combining
\cref{eq:lrm-high-block-maxima,eq:lrm-low-block-maxima},
\[
 B_\ell:=
 \sup_{\tau_\ell<k\leq\tau_{\ell+1}}
 \left\|\sum_{n=\tau_\ell+1}^{k}
 \frac{\varepsilon_n}{\sqrt{S_{n-1}}}\right\|
 \leq B_{\ell,H}+B_{\ell,L}\longrightarrow0
 \quad\text{a.s.}
\]

\paragraph{Step 4: arbitrary windows and the boundary increment.}
Set $e_n:=\|\varepsilon_n\|/\sqrt{S_{n-1}}$. On $\mathcal A^c$,
objective stability and \cref{lem:nonflat-linear-growth} give bounded
iterates. Global smoothness therefore gives
$\sup_n\|h_n\|<\infty$. By \cref{eq:lrm-relative-jump},
\begin{equation}
\label{eq:lrm-boundary-increment}
    e_n
    \leq
    r_n+\frac{\|h_n\|}{\sqrt{S_{n-1}}}
    \longrightarrow0
    \qquad\text{a.s. on }\mathcal A^c.
\end{equation}
For $n\geq1$, let $\ell(n):=\lfloor t_n/T\rfloor$. Then
$\tau_{\ell(n)}\leq n<\tau_{\ell(n)+1}$. If
$k\leq m(t_n+T)$, then $k<\tau_{\ell(n)+2}$. A sum beginning strictly
inside the $\ell(n)$-th block is the difference of two block prefixes; if
$n=\tau_{\ell(n)}$, its first summand belongs to the preceding block and is
bounded by $e_n$. Therefore,
\begin{equation}
\label{eq:lrm-arbitrary-window}
    \sup_{n\leq k\leq m(t_n+T)}
    \left\|
      \sum_{i=n}^{k}\frac{\varepsilon_i}{\sqrt{S_{i-1}}}
    \right\|
    \leq
    e_n+2B_{\ell(n)}+B_{\ell(n)+1}.
\end{equation}
Because $t_n\to\infty$, one has $\ell(n)\to\infty$.  Hence
\cref{eq:lrm-boundary-increment,eq:lrm-arbitrary-window} show that the
predictable-denominator martingale vanishes on every finite effective-time
window. Combining this with the absolutely summable correction in
\cref{eq:lrm-denominator-decomposition,eq:lrm-summable-correction} and
multiplying by $\alpha_0$ proves \cref{eq:lrm-finite-window} for the fixed
$T$. Intersecting the probability-one events for positive rational $T$ and
bounding an arbitrary real $T>0$ by a larger rational window gives the
conclusion simultaneously for every $T>0$.
\end{proof}

\section{Proof of Asymptotic Stationarity}
\label{appendix:asymptotic}

\begin{namedproof}{Proof of Theorem~\ref{thm:main}, part~\ref{thm:main:stationarity}.}
Let \(\mathcal A:=\{\lim_{n\to\infty}S_n<\infty\}\).

\paragraph{Case 1: \(\mathcal A\).}
By \cref{eq:gradient-energy-standard}, for every \(\delta>0\),
\[
    \sum_{n=1}^{\infty}
    \ind_{\{\|\nabla g(\theta_n)\|>\delta\}}
    \frac{\|\nabla g(\theta_n)\|^2}{S_{n-1}}
    <
    \infty
    \qquad\text{almost surely}.
\]
On \(\mathcal A\), \(S_{n-1}\leq S_\infty<\infty\). Hence, whenever
\(\|\nabla g(\theta_n)\|>\delta\),
\[
    \frac{\|\nabla g(\theta_n)\|^2}{S_{n-1}}
    \geq
    \frac{\delta^2}{S_\infty}.
\]
Therefore, the event
\(\{\|\nabla g(\theta_n)\|>\delta\}\) can occur only finitely often.
Since this holds for every rational \(\delta>0\),
\[
    \|\nabla g(\theta_n)\|\longrightarrow0
    \qquad\text{almost surely on }\mathcal A.
\]

\paragraph{Case 2: \(\mathcal A^c\).}
Here \(S_n\to\infty\). Recall \cref{eq:oracle-noise} and set
\(\gamma_n:=\alpha_0/\sqrt{S_n}\). Then
\[
    \theta_{n+1}
    =
    \theta_n
    -
    \gamma_n\bigl(\nabla g(\theta_n)+\varepsilon_n\bigr).
\]
By \cref{lem:effective-time-divergence}, $    \gamma_n\to0,
\sum_{n=1}^{\infty}\gamma_n=\infty.$
Iterate boundedness follows from
\cref{thm:main}~\ref{thm:main:boundedness}, while
\cref{inequ:gammaU} gives the finite-window perturbation condition.
Consequently, the continuous-time interpolation is an asymptotic
pseudotrajectory of \(\dot z=-\nabla g(z)\).

The objective \(g\) is a strict Lyapunov function for this flow. The ODE
limit-set principle together with the weak Sard condition implies that every
accumulation point of
\(\{\theta_n\}\) is stationary. Since the sequence is bounded, its
limit set is nonempty and compact. If
\(\|\nabla g(\theta_n)\|\not\to0\), there would exist
\(\varepsilon>0\) and a convergent subsequence
\(\theta_{n_j}\to\bar{\theta}\) such that
\(\|\nabla g(\theta_{n_j})\|\geq\varepsilon\). Continuity would then
give \(\|\nabla g(\bar{\theta})\|\geq\varepsilon\), contradicting
stationarity of every accumulation point. Hence,
\[
    \|\nabla g(\theta_n)\|\longrightarrow0
    \qquad\text{almost surely on }\mathcal A^c.
\]
Combining the two cases proves the almost-sure conclusion.

For the mean-square conclusion, \cref{eq:smooth-gradient-objective-bound} and
\cref{thm:main}~\ref{thm:main:stability} give
\(\sup_{n\geq1}\|\nabla g(\theta_n)\|^2
\leq2L\sup_{n\geq1}g(\theta_n)\), and the right-hand side is integrable.
Dominated convergence therefore yields $ \lim_{n\to\infty}\E\|\nabla g(\theta_n)\|^2=0.$
\end{namedproof}

\section{Proofs for Power-Normalized AdaGrad with \(q>1/2\)}
\label{app:power-normalization-proofs}

We prove \cref{thm:power-affine,thm:power-uniform-l2}.

Throughout this section, consider \cref{eq:power-normalized-family} and write
\[
\begin{gathered}
 h_n:=\nabla g(\theta_n),\qquad P_n:=S_{n-1},\qquad
 \Delta_n:=\|G_n\|^2,\qquad Q_n:=S_n=P_n+\Delta_n,\\
 a_n:=P_n^{-q},\qquad b_n:=Q_n^{-q},\qquad
 \chi_n:=1-(P_n/Q_n)^q,\qquad d_n:=a_n-b_n=a_n\chi_n,\\
 \zeta_n:=a_n\|h_n\|^2,\qquad
 \kappa:=\frac{\alpha_0\beta}{2},\qquad
 V_n:=g(\theta_n)+\kappa\zeta_n.
\end{gathered}
\]
Then \(P_{n+1}=Q_n\), \(\zeta_{n+1}=b_n\|h_{n+1}\|^2\), and
\(P_n,h_n,a_n,\zeta_n,V_n\) are \(\mathscr F_{n-1}\)-measurable.

\subsection{A Conditional Lyapunov Inequality}

\begin{lem}[One-step power-normalized Lyapunov inequality]
\label{lem:power-one-step}
Under \cref{eq:power-unbiased,eq:power-affine-moment}, there are constants
$c_p>0$ and $C_p<\infty$, depending only on
$q,L,\alpha_0,\beta,\sigma^2$, such that
\begin{equation}
\label{eq:power-one-step}
    \E[V_{n+1}\mid\mathscr F_{n-1}]
    \leq
    V_n-c_p\zeta_n
    +C_p\E\left[
       \frac{\Delta_n}{Q_n^{2q}}
       +\frac{\Delta_n}{Q_n^{3q}}
       +d_n
       \;\middle|\;\mathscr F_{n-1}
    \right].
\end{equation}
One may take $c_p=\alpha_0/4$ after enlarging $C_p$.
\end{lem}

\begin{proof}
The descent lemma gives
\begin{equation}
\label{eq:power-descent}
    g(\theta_{n+1})-g(\theta_n)
    \leq
    -\alpha_0b_n\langle h_n,G_n\rangle
    +\frac{L\alpha_0^2}{2}\frac{\Delta_n}{Q_n^{2q}}.
\end{equation}
Since $b_n=a_n(1-\chi_n)$, conditional unbiasedness yields
\[
\E[-b_n\langle h_n,G_n\rangle\mid\mathscr F_{n-1}]
=-a_n\|h_n\|^2
+a_n\E[\chi_n\langle h_n,G_n\rangle\mid\mathscr F_{n-1}].
\]
Conditional Cauchy--Schwarz, Young's inequality, and
\cref{eq:power-affine-moment} imply
\[
\begin{aligned}
&a_n\E[\chi_n\langle h_n,G_n\rangle\mid\mathscr F_{n-1}]\leq
\frac12\zeta_n
+\frac{\beta}{2}\zeta_n
  \E[\chi_n^2\mid\mathscr F_{n-1}]
+\frac{\sigma^2}{2}a_n
  \E[\chi_n^2\mid\mathscr F_{n-1}].
\end{aligned}
\]
Because $0\leq\chi_n\leq1$,
$a_n\chi_n^2\leq a_n\chi_n=d_n$. Adding
$\kappa(\E[\zeta_{n+1}\mid\mathscr F_{n-1}]-\zeta_n)$ to
\cref{eq:power-descent} and using $\kappa=\alpha_0\beta/2$ gives
\begin{align}
\E[V_{n+1}\mid\mathscr F_{n-1}]-V_n
\leq{}&
-\frac{\alpha_0}{2}\zeta_n
+\kappa\E\left[
  \bigl(\|h_{n+1}\|^2-\|h_n\|^2\bigr)b_n
  \;\middle|\;\mathscr F_{n-1}
\right]
\notag\\
&+
\frac{\alpha_0\sigma^2}{2}
\E[d_n\mid\mathscr F_{n-1}]
+
\frac{L\alpha_0^2}{2}
\E\left[\frac{\Delta_n}{Q_n^{2q}}\middle|\mathscr F_{n-1}\right].
\label{eq:power-intermediate}
\end{align}
Global smoothness gives
\(\|h_{n+1}-h_n\|\leq
L\|\theta_{n+1}-\theta_n\|=L\alpha_0\sqrt{\Delta_n}\,Q_n^{-q}\).
Hence
\[
\bigl(\|h_{n+1}\|^2-\|h_n\|^2\bigr)b_n
\leq
2L\alpha_0\|h_n\|\sqrt{\Delta_n}\,Q_n^{-2q}
+L^2\alpha_0^2\frac{\Delta_n}{Q_n^{3q}}.
\]
For every $\varepsilon>0$, Young's inequality and $P_n\leq Q_n$ give
\[
2L\alpha_0\|h_n\|\sqrt{\Delta_n}\,Q_n^{-2q}
\leq
\varepsilon\zeta_n
+
\frac{L^2\alpha_0^2}{\varepsilon}
\frac{\Delta_n}{Q_n^{3q}}.
\]
If $\beta>0$, choose $\varepsilon=(2\beta)^{-1}$ so that
$\kappa\varepsilon=\alpha_0/4$; if $\beta=0$, the term multiplied by
$\kappa$ vanishes. Substitution into \cref{eq:power-intermediate} proves
\cref{eq:power-one-step}.
\end{proof}

For later use, if \(r>1\), monotonicity of \(x^{-r}\) and
\(P_{n+1}=Q_n\) give, for every \(N\),
\begin{equation}
\label{eq:power-budget-general}
 \sum_{n=1}^{N}\frac{\Delta_n}{Q_n^r}
 \leq\int_{P_1}^{Q_N}x^{-r}\,\mathrm dx
 \leq\frac{P_1^{1-r}}{r-1}.
\end{equation}
Moreover, \(d_n=P_n^{-q}-P_{n+1}^{-q}\), so
\begin{equation}
\label{eq:power-d-budget}
 \sum_{n=1}^{N}d_n=P_1^{-q}-Q_N^{-q}\leq P_1^{-q}.
\end{equation}

\subsection{Stability in the Summable Regime}

\begin{namedproof}{Proof of the stability assertions in \cref{thm:power-affine}.}
Set
\[
    R_n
    :=
    \frac{\Delta_n}{Q_n^{2q}}
    +\frac{\Delta_n}{Q_n^{3q}}
    +d_n.
\]
Since $q>1/2$, \cref{eq:power-budget-general,eq:power-d-budget} give the pathwise bound
\[
    \sum_{n=1}^{\infty}R_n
    \leq
    \frac{P_1^{1-2q}}{2q-1}
    +\frac{P_1^{1-3q}}{3q-1}
    +P_1^{-q}
    =:C_R<\infty.
\]
Tonelli's theorem yields
\[
    \E\left[
      \sum_{n=1}^{\infty}
      \E[R_n\mid\mathscr F_{n-1}]
    \right]
    =\E\left[\sum_{n=1}^{\infty}R_n\right]
    \leq C_R.
\]
In particular, the predictable error series is finite almost surely.
To justify the almost-supermartingale argument, $V_n$ is adapted to the shifted filtration
$\mathcal G_n:=\mathscr F_{n-1}$ and is integrable at each finite time.
Indeed, for $q>1/2$ the quantity
$u^{1/2}(P_1+u)^{-q}$ is bounded on $[0,\infty)$, so every update length is
deterministically bounded; global Lipschitz continuity of $\nabla g$ then
bounds $V_n$ deterministically at each fixed index. Applying the
Robbins--Siegmund almost-supermartingale theorem
\citep{robbins1971convergence} to \cref{eq:power-one-step} shows that
$V_n$ converges to a finite random variable and
\begin{equation}
\label{eq:power-zeta-sum}
    \sum_{n=1}^{\infty}\zeta_n<\infty
    \qquad\text{almost surely}.
\end{equation}
Hence $\zeta_n\to0$ and
$g(\theta_n)=V_n-\kappa\zeta_n$ converges to a finite random variable.
Since $\zeta_n=S_{n-1}^{-q}\|\nabla g(\theta_n)\|^2$,
\cref{eq:power-as-stability} follows.
\end{namedproof}

\begin{namedproof}{Proof of \cref{eq:power-integrable-stability} in \cref{thm:power-uniform-l2}.}
Under \cref{eq:power-uniform-second-moment}, conditional Jensen gives
$\|h_n\|^2\leq K$. We may therefore take $\beta=0$ and
$\sigma^2=K$, so no gradient-dependent correction is needed. Define
\begin{equation}
\label{eq:power-uniform-martingale}
    M_n
    :=
    \alpha_0\left(
      \E[b_n\langle h_n,G_n\rangle\mid\mathscr F_{n-1}]
      -b_n\langle h_n,G_n\rangle
    \right).
\end{equation}
Repeating the preceding descent calculation with $\beta=0$ gives the
pathwise bound
\begin{equation}
\label{eq:power-uniform-pathwise-descent}
    g(\theta_{n+1})-g(\theta_n)
    \leq
    -\frac{\alpha_0}{2}\zeta_n
    +\frac{\alpha_0K}{2}\E[d_n\mid\mathscr F_{n-1}]
    +\frac{L\alpha_0^2}{2}\frac{\Delta_n}{Q_n^{2q}}
    +M_n.
\end{equation}
Moreover,
\[
    \E[M_n^2\mid\mathscr F_{n-1}]
    \leq
    \alpha_0^2\|h_n\|^2
    \E\left[\frac{\Delta_n}{Q_n^{2q}}\middle|\mathscr F_{n-1}\right]
    \leq
    \alpha_0^2K
    \E\left[\frac{\Delta_n}{Q_n^{2q}}\middle|\mathscr F_{n-1}\right].
\]
Let $\mathcal M_N:=\sum_{n=1}^{N}M_n$. By
\cref{eq:power-budget-general} and Tonelli's theorem,
\[
    \sup_{N\geq1}\E[\mathcal M_N^2]
    \leq
    \frac{\alpha_0^2KP_1^{1-2q}}{2q-1}
    =:C_M<\infty.
\]
Thus $\mathcal M_N$ is $L^2$-bounded. Doob's $L^2$ maximal inequality and
Cauchy--Schwarz give
\begin{equation}
\label{eq:power-uniform-martingale-max}
    \E\left[\sup_{N\geq1}|\mathcal M_N|\right]
    \leq 2\sqrt{C_M}<\infty.
\end{equation}
Summing \cref{eq:power-uniform-pathwise-descent} through $N$ and dropping the
negative drift gives
\[
\begin{aligned}
    \sup_{N\geq1}g(\theta_{N+1})
    \leq{}&g(\theta_1)
    +\frac{\alpha_0K}{2}
      \sum_{n=1}^{\infty}\E[d_n\mid\mathscr F_{n-1}]+\frac{L\alpha_0^2}{2}
      \sum_{n=1}^{\infty}\frac{\Delta_n}{Q_n^{2q}}
    +\sup_{N\geq1}|\mathcal M_N|.
\end{aligned}
\]
The first two nonnegative series on the right have finite expectations by
\cref{eq:power-d-budget,eq:power-budget-general} and Tonelli's theorem; the
last term is integrable by \cref{eq:power-uniform-martingale-max}.  This proves
\cref{eq:power-integrable-stability}.
\end{namedproof}

\subsection{Effective Time and Current-Centered Noise}

\begin{lem}[Divergence of the predictable effective time]
\label{lem:power-effective-time}
Assume the hypotheses of \cref{thm:power-affine} and let
$q\in(1/2,1]$. On the event $\{P_n\to\infty\}$,
\begin{equation}
\label{eq:power-effective-time}
    \sum_{n=1}^{\infty}a_n
    =
    \sum_{n=1}^{\infty}S_{n-1}^{-q}
    =\infty
    \qquad\text{almost surely}.
\end{equation}
\end{lem}

\begin{proof}
The objective bound from \cref{eq:power-as-stability} and
\(\|h_n\|^2\leq2Lg(\theta_n)\) imply
\(\sup_n\|h_n\|<\infty\) almost surely. Hence, pathwise,
\(m_n:=\E[\Delta_n\mid\mathscr F_{n-1}]
\leq\beta\sup_j\|h_j\|^2+\sigma^2=:K(\omega)<\infty\).
Suppose, for contradiction, that $\sum_na_n<\infty$. Then
$A_\infty:=\sum_na_nm_n<\infty$. By the conditional summability lemma
\cref{lem:conditional-summability}, applied after localization if necessary,
\begin{equation}
\label{eq:power-an-delta-finite}
    \sum_{n=1}^{\infty}a_n\Delta_n<\infty
    \qquad\text{almost surely}.
\end{equation}
On the other hand, since \(P_{n+1}=P_n+\Delta_n\) and
\(x\mapsto x^{-q}\) is decreasing,
\(a_n\Delta_n\geq\int_{P_n}^{P_{n+1}}x^{-q}\,\mathrm dx\).
For $q<1$ the cumulative integral grows as $P_{N+1}^{1-q}/(1-q)$, and
for $q=1$ it grows as $\log P_{N+1}$. Both diverge on
$\{P_n\to\infty\}$, contradicting \cref{eq:power-an-delta-finite}.
\end{proof}

\begin{lem}[Current-centered finite-window noise]
\label{lem:power-current-centered}
Under the hypotheses of \cref{lem:power-effective-time}, define
\begin{equation}
\label{eq:power-Z-r}
    Z_n:=b_nG_n-\E[b_nG_n\mid\mathscr F_{n-1}],
    \qquad
    r_n:=\E[b_nG_n\mid\mathscr F_{n-1}]-a_nh_n.
\end{equation}
On $\{P_n\to\infty\}$, the vector martingale
$\sum_{j=1}^{n}Z_j$ converges almost surely and
\begin{equation}
\label{eq:power-r-small}
    \frac{\|r_n\|}{a_n}\longrightarrow0
    \qquad\text{almost surely}.
\end{equation}
Consequently, for every $T>0$, if $ k(n,T)
    :=
    \max\left\{k\geq n:
      \alpha_0\sum_{j=n}^{k}a_j\leq T
    \right\},$
then
\begin{equation}
\label{eq:power-finite-window}
    \lim_{n\to\infty}
    \sup_{n\leq k\leq k(n,T)}
    \left\|
      \alpha_0\sum_{j=n}^{k}(Z_j+r_j)
    \right\|
    =0
    \qquad\text{almost surely}.
\end{equation}
The same conclusion holds under the usual ODE convention that includes the
single increment crossing the right endpoint.
\end{lem}

\begin{proof}
The process \(Z_n\) is an \(\R^d\)-valued martingale difference and
\(\E[\|Z_n\|^2\mid\mathscr F_{n-1}]
\leq\E[\Delta_n/Q_n^{2q}\mid\mathscr F_{n-1}]\).
By \cref{eq:power-budget-general} and Tonelli's theorem, the total expected
quadratic variation is finite. Hence the martingale
$\sum_{j=1}^{n}Z_j$ is $L^2$-bounded and converges almost surely.

Next, since $d_n=a_n\chi_n$ and $q\leq1$,
\[
    \|r_n\|
    \leq
    a_n\E[\chi_n\sqrt{\Delta_n}\mid\mathscr F_{n-1}],
    \qquad
    \chi_n
    =1-\left(\frac{P_n}{Q_n}\right)^q
    \leq\frac{\Delta_n}{Q_n}.
\]
Moreover,
\[
    \chi_n\sqrt{\Delta_n}
    \leq
    \frac{\Delta_n^{3/2}}{Q_n}
    \leq
    \frac{\Delta_n}{\sqrt{P_n}},
\]
so, with the pathwise bound $m_n\leq K(\omega)$ from the proof of
\cref{lem:power-effective-time},
\[
    \frac{\|r_n\|}{a_n}
    \leq\frac{m_n}{\sqrt{P_n}}
    \leq\frac{K(\omega)}{\sqrt{P_n}}\to0.
\]
This proves \cref{eq:power-r-small}.  The martingale tail is uniformly Cauchy,
whereas the bias contribution over the above window is bounded by
\[
    \alpha_0\sum_{j=n}^{k(n,T)}\|r_j\|
    \leq
    T\sup_{j\geq n}\frac{\|r_j\|}{a_j}\to0,
\]
which proves \cref{eq:power-finite-window} for windows that do not include
the crossing increment.

It remains to control that single boundary increment. Since
$Z_n+r_n=b_nG_n-a_nh_n$, write $u=\Delta_n/P_n$. Then
\[
    b_n\|G_n\|
    =
    P_n^{1/2-q}\frac{\sqrt u}{(1+u)^q}
    \leq C_qP_n^{1/2-q}\to0
\]
because $q>1/2$. The objective bound makes $\|h_n\|$ bounded, so
$a_n\|h_n\|\to0$ as well. Adding or deleting the crossing increment therefore
does not change the limit.
\end{proof}

\subsection{Stationarity}

\begin{namedproof}{Proof of \cref{eq:power-as-stationarity}.}
Assume $q\in(1/2,1]$. On $\{\sup_nP_n<\infty\}$,
\cref{eq:power-zeta-sum} and
$a_n\geq(\sup_jP_j)^{-q}>0$ imply
$\sum_n\|h_n\|^2<\infty$, hence $\|h_n\|\to0$.

On $\{P_n\to\infty\}$, \cref{eq:power-Z-r} rewrites the recursion as
\begin{equation}
\label{eq:power-sa-recursion}
    \theta_{n+1}
    =
    \theta_n-\alpha_0a_nh_n-\alpha_0(Z_n+r_n).
\end{equation}
By \cref{lem:power-effective-time}, $\sum_na_n=\infty$ and $a_n\to0$;
\cref{lem:power-current-centered} supplies the finite-window perturbation
condition. We use these facts directly; no precompactness assumption on the iterates is
needed for this argument.

First,
\begin{equation}
\label{eq:power-liminf-gradient}
    \liminf_{n\to\infty}\|h_n\|=0,
\end{equation}
because otherwise the series
\(\sum_na_n\|h_n\|^2\) in \cref{eq:power-zeta-sum} would diverge. Also,
\(\|h_{n+1}-h_n\|\leq L\alpha_0b_n\|G_n\|\to0\),
where the last limit was established in the boundary-increment argument in
\cref{lem:power-current-centered}.

Suppose that $\|h_n\|$ does not converge to zero. Choose $\varepsilon>0$
such that $\limsup_n\|h_n\|>3\varepsilon$. Using
\cref{eq:power-liminf-gradient} and the vanishing one-step gradient change,
construct pairwise disjoint intervals $[s_j,t_j]$ with $s_j\to\infty$ such
that, for all sufficiently large $j$,
\begin{equation}
\label{eq:power-upcrossing-gradient-band}
    \varepsilon<\|h_n\|<3\varepsilon
    \quad(s_j\leq n<t_j),
    \qquad
    \|h_{s_j}\|\leq\frac32\varepsilon,
    \qquad
    \|h_{t_j}\|\geq3\varepsilon.
\end{equation}
Set
\[
    A_j:=\alpha_0\sum_{n=s_j}^{t_j-1}a_n,
    \qquad
    T_0:=\frac1{6L}.
\]
If $A_j\leq T_0$ along an infinite subsequence, then
\cref{eq:power-finite-window,eq:power-sa-recursion,eq:power-upcrossing-gradient-band}
give
\[
\begin{aligned}
    \frac32\varepsilon
    &\leq\|h_{t_j}-h_{s_j}\|
    \leq L\|\theta_{t_j}-\theta_{s_j}\|\leq
    L\left(
       3\varepsilon A_j
       +\alpha_0\left\|\sum_{n=s_j}^{t_j-1}(Z_n+r_n)\right\|
    \right)
    \leq\frac{\varepsilon}{2}+o(1),
\end{aligned}
\]
a contradiction. Thus $A_j>T_0$ eventually. Since the intervals are
disjoint, \cref{eq:power-upcrossing-gradient-band} then gives
\[
    \sum_{n=s_j}^{t_j-1}a_n\|h_n\|^2
    \geq
    \frac{\varepsilon^2T_0}{\alpha_0}
\]
for infinitely many $j$, contradicting \cref{eq:power-zeta-sum}. Therefore
$\|\nabla g(\theta_n)\|\to0$ almost surely, proving
\cref{eq:power-as-stationarity}.
\end{namedproof}

\begin{namedproof}{Proof of \cref{eq:power-ms-stationarity}.}
Under \cref{eq:power-uniform-second-moment}, conditional Jensen gives the
deterministic domination $\|\nabla g(\theta_n)\|^2\leq K$. For $q\leq1$,
\cref{eq:power-as-stationarity} gives almost-sure convergence to zero.
Dominated convergence proves \cref{eq:power-ms-stationarity}.
\end{namedproof}

\bibliography{adagrad_ref_master}

\end{document}